\documentclass[12pt]{article}
\usepackage{amsfonts}
\usepackage{stmaryrd}
\usepackage{stmaryrd}
\usepackage{graphicx}
\usepackage{graphics}
\usepackage{epstopdf}
\usepackage{epsfig}
\usepackage{amsmath}
\usepackage{amsthm}
\usepackage{amssymb}
\usepackage{float}
\usepackage{cite}
\usepackage{multirow}
\usepackage{verbatim}
\usepackage{fancyhdr}
\usepackage{subfigure}
\usepackage{color}
\usepackage{mathtools}
\usepackage{sectsty}
\usepackage[title]{appendix}
\usepackage{threeparttable}
\usepackage{dcolumn}
\usepackage{booktabs}
\usepackage{indentfirst}
\usepackage{setspace}
\usepackage{bm}
\usepackage{enumerate}
\usepackage[labelfont=bf,labelsep=period]{caption}
\usepackage{setspace}
\usepackage{diagbox}
\usepackage{hyperref}
\hypersetup{colorlinks,linkcolor=blue,citecolor=blue}

\newcolumntype{d}[1]{D{.}{.}{#1}}
\theoremstyle{definition}

\numberwithin{equation}{section}

\makeatletter
\renewcommand{\maketag@@@}[1]{\hbox{\m@th\normalsize\normalfont#1}}
\makeatother

\allowdisplaybreaks
\begin{document}
\pagenumbering{arabic}
\baselineskip=1.4pc

\vspace*{0.5in}

\begin{center}

{{\bf \Large
Inverse Lax-Wendroff boundary treatment for solving conservation laws with finite difference HWENO methods}}

\end{center}

\vspace{.03in}
\vspace{.03in}

\centerline{
Guanghao Zhu\footnote{School of Mathematical Sciences, 
    University of Science and Technology of China, 
    Hefei, Anhui 230026, China. 
    E-mail: zhugy@mail.ustc.edu.cn.},
Yan Jiang\footnote{School of Mathematical Sciences, 
    University of Science and Technology of China, 
    Hefei, Anhui 230026, China. 
    E-mail: jiangy@ustc.edu.cn. 
    Research is supported in part by National Natural Science Foundation of China grant 12271499.}
    \footnote{Laoshan Laboratory, Qingdao 266237, China},
Zhuang Zhao\footnote{School of Mathematical Sciences and Fujian Provincial Key Laboratory of Mathematical Modeling and High-Performance Scientific Computing, 
    Xiamen University, Xiamen, Fujian 361005, China. 
    E-mail: zzhao@xmu.edu.cn. 
    Research is supported in part by National Natural Science Foundation of Xiamen, China grant 3502Z202472004 and National Natural Science Foundation of China grant 12401541.} 
and Mengping Zhang\footnote{School of Mathematical Sciences,
    University of Science and Technology of China, Hefei,
    Anhui 230026, China.  
    E-mail: mpzhang@ustc.edu.cn.}
}

\vspace{.1in}

\noindent
\textbf{Abstract:}
This paper presents a novel inverse Lax-Wendroff (ILW) boundary treatment for finite difference Hermite weighted essentially non-oscillatory (HWENO) schemes to solve hyperbolic conservation laws on arbitrary geometries.
The complex geometric domain is divided by a uniform Cartesian grid, resulting in challenge in boundary treatment. The proposed ILW boundary treatment could provide high order approximations of both solution values and spatial derivatives at ghost points outside the computational domain.
Distinct from existing ILW approaches, our boundary treatment constructs the extrapolation via optimized through a least squares formulation, coupled with the spatial derivatives at the boundary obtained via the ILW procedure. 
Theoretical analysis indicates that compared with other ILW methods, our proposed one would require fewer terms by using the relatively complicated ILW procedure and thus improve computational efficiency while preserving accuracy and stability.
The effectiveness and robustness of the method are validated through numerical experiments.

\bigskip


\pagenumbering{arabic}

\baselineskip=1.4pc

\section{Introduction}

Hyperbolic conservation laws pose significant challenges in numerical simulations, particularly when complex geometries require accurate and stable boundary treatments. The high order finite difference Hermite weighted essentially non-oscillatory (HWENO) schemes have emerged as powerful tools for resolving discontinuities and sharp gradients. 
However, their application to complex geometric domains remains hindered by unresolved computational challenges.
In this paper, we propose an efficient boundary treatment for the fifth order finite difference HWENO scheme \cite{fan2023robust} to solve hyperbolic conservation laws, specifically tailored to address these limitations in irregular domains.

Distinct from traditional weighted essentially non-oscillatory (WENO) schemes, the HWENO scheme introduces derivative equations into its formulation. This strategy utilizes not only the solution values but also their first derivatives, enabling more compact stencils, enhancing resolution, and reducing computational costs at equivalent orders of accuracy compared to the WENO schemes using same nonlinear weight definitions. 
The HWENO scheme was first proposed by Qiu and Shu \cite{qiu2004hermite, qiu2005hermite} in the finite volume framework as a limiter for the discontinuous Galerkin (DG) method. 
Two high order finite difference HWENO schemes were later developed in \cite{liu2015finite,liu2016finite}. However, due to the presence of mixed derivatives, both finite difference methods reduce to fourth order accuracy in two-dimensional cases. Subsequently, Zhao et al. \cite{zhao2020modified} took the idea of limiter from the DG method to modify the  first-order derivatives by applying limiters, enabling the algorithm to achieve fifth-order accuracy in two-dimensional cases. On this basis, the HWENO-L scheme was further proposed in \cite{zhang2023fifth}, which reduced computational storage and optimized time discretization to improve compactness. Recently, the HWENO-R scheme by Fan et al. \cite{fan2023robust} represented  a significant advancement through optimized reconstruction stencils, demonstrating superior robustness for extreme flow conditions.
For more details about the HWENO scheme, please refer to \cite{zahran2016seventh, du2018hermite, li2021multi, Fan2024OEHWENO, xie2025113673}. 
However, to the best of our knowledge, there is currently no systematic research on the numerical boundary treatments specifically for the HWENO scheme, particularly for complex geometries.
This gap motivates our work, which focuses on extending the applicability of the finite difference HWENO scheme to complex geometries through proper numerical boundary treatments.

Within the finite difference framework, two principal challenges arise when implementing HWENO schemes for complex boundaries. 
First, the HWENO scheme depends on information from neighboring points, which necessitates the construction of values and first derivatives at ghost points near the boundary. Numerical experiments have demonstrated that improper handling of these ghost points can result in loss of accuracy and instability. 
Secondly, grid points generally do not coincide with physical boundaries, requiring specialized techniques to preserve accuracy while imposing the boundary conditions. 
Existing boundary treatment strategies broadly fall into two categories: body-fitted mesh methods and non-body-fitted mesh methods. 
The body-fitted mesh methods ensure that the grid points fall on the boundary, and the numerical solution can accurately satisfy the boundary conditions and achieve high accuracy. However, the mesh  generation requires additional computational cost, particularly for problems involving moving boundaries. 
The non-body-fitted mesh methods can use uniform Cartesian grids, with the advantages of simple data structure and ease of computation. 
However, the misalignment between the boundary and grid points often results in the so called ``cut-cell'' problem, which imposes stricter constraints on time step sizes to ensure numerical stability.
To overcome these challenges, a lot of methods have been proposed, including the immersed boundary method \cite{peskin1972flow} and the inverse Lax-Wendroff (ILW) method. 
It should be noted that the immersed boundary method typically necessitates additional treatment to attain high-order accuracy.
In this work, we adopt the ILW method, which naturally achieves high-order accuracy and effectively avoids the cut-cell problem.

The ILW method was proposed by Tan and Shu \cite{tan2010inverse} as a high-order boundary treatment method for solving hyperbolic conservation laws with finite difference WENO schemes, and it has been applied to the simulation of compressible flows in complex domains. This ILW procedure transforms the normal derivatives on the boundary into time derivatives and tangential derivatives utilizing the governing equations repeatedly. With these normal derivatives, the accurate values of ghost points near the boundary are imposed by a Taylor expansion. This method can achieve arbitrary higher-order accuracy at the boundary. 
To avoid complex calculations caused by the ILW procedure for high-order derivatives, Tan et al.\cite{tan2012efficient} introduced the simplified ILW (SILW), which employs extrapolation to obtain those high-order normal derivatives. A stability analysis for SILW was provided in \cite{li2016stability}, demonstrating that this method effectively avoids the cut-cell problem and specifying the number of low-order terms required via the ILW procedure.  
Recently, Liu et al. \cite{liu2024113259} introduced a new SILW method for the finite different method to further simplify calculations and improve computational efficiency, in which the number of low-order terms required with the ILW procedure is reduced while high-order accuracy and stability are maintained. 
The SILW method was also extended to the steady state of hyperbolic conservation laws \cite{li2023fixed}, the convection-diffusion equations \cite{lu2016inverse, li2017stability, li2022stability}, the Boltzmann equation\cite{filbet2013inverse} and the moving boundary problems \cite{cheng2021high, liu2022numerical, liu2023high, LIU2025113942}. All these works are based on the finite difference methods. So far, the ILW boundary treatment has been designed to solve the conservation laws with finite volume method \cite{zhu2024inverse} and discontinuous Galerkin method \cite{yang2024inverse} as well.  

In this paper, we will propose a novel integration of the ILW method into the fifth order finite different HWENO-R scheme to address high-order ghost point reconstruction for both solution values and derivatives, ensuring stability in irregular domains. 
This framework will fill a critical methodological void in HWENO applications on complex geometric domains. 
In particular, a new approximation based on the least square method instead of interpolation is designed in order to use less terms with the complex ILW procedure and reduce computational complexity.
Stability analysis will be provided to validate the robustness of the proposed algorithm. 
Without specified clarification, the notation ``HWENO'' in the following refers to the HWENO-R scheme.

The organization of this paper is as follows. In Section 2, we review the framework of the finite difference HWENO scheme. In Section 3, we will describe the ILW method for one-dimensional scalar conservation laws and perform a linear stability analysis using the characteristic spectral visualization method. The algorithm will be extended to two-dimensional systems in Section 4. Section 5 presents numerical experiments to validate the efficiency and high-order accuracy of the proposed algorithm. Finally, conclusions are provided in Section 6.

\section{Finite difference HWENO scheme}
In this section, we provide a brief review on the finite difference HWENO scheme for hyperbolic conservation laws \cite{fan2023robust}. 
We first consider the one-dimensional scalar case:
\begin{equation}\label{ODSE-1}
    u_t+f(u)_x=0.
\end{equation}
In the finite difference framework, the grid is uniformly divided into mesh size $\Delta x=x_{i+1}-x_{i}$. 
Defining $v=u_x$ and $h(u,v)=f(u)_x=f'(u)u_x=f'(u)v$, we then design the HWENO scheme by adding an additional derivative equation for (\ref{ODSE-1}):
\begin{equation}
	\begin{cases}
		u_t+f(u)_x=0, \\
		v_t+h(u,v)_x=0.
	\end{cases}
\end{equation}
The conservative semi-discrete HWENO scheme can be written as:
\begin{equation}
\begin{cases}
    \frac{d}{dt}u_i(t)=-\frac{1}{\Delta x}(\hat{f}_{i+\frac{1}{2}}-\hat{f}_{i-\frac{1}{2}}), \\
    \frac{d}{dt}v_i(t)=-\frac{1}{\Delta x}(\hat{h}_{i+\frac{1}{2}}-\hat{h}_{i-\frac{1}{2}}).
\end{cases}
\label{ODSE2}
\end{equation}
Here, $u_i(t)$ and $v_i(t)$ are the approximation of point values $u(x_i, t)$ and $v(x_i,t)$, respectively. $\hat{f}_{i+\frac{1}{2}}$ and $\hat{h}_{i+\frac{1}{2}}$ are the numerical fluxes at half points, which can be obtained via the HWENO reconstruction. To achieve a $d$-th order accuracy with respect to $u$, we require 
\begin{equation*}
\begin{aligned}
    & \mathcal{F}_i(u,v) := \frac{1}{\Delta x}(\hat{f}_{i+\frac{1}{2}}-\hat{f}_{i-\frac{1}{2}}) = f(u)_x|_{x_i} + \mathcal{O}(\Delta x^d), \\
    & \mathcal{H}_i(u,v) := \frac{1}{\Delta x}(\hat{h}_{i+\frac{1}{2}}-\hat{h}_{i-\frac{1}{2}}) = h(u,v)_x|_{x_i} + \mathcal{O}(\Delta x^{d-1}).
\end{aligned}
\end{equation*}

\noindent
Moreover, to ensure the stability, it is important to use upwinding strategy when constructing the fluxes. One common practice is to use flux splitting. We do a flux splitting on $f$ and $h$ respectively,
	\begin{equation*}
		\begin{aligned}
		f^{\pm}_i=f^{\pm}(u_i)=\frac{1}{2}(f(u_i)\pm\alpha u_i),\quad 
		h^{\pm}_i=h^{\pm}(u_i,v_i)=\frac{1}{2}(h(u_i,v_i)\pm\alpha v_i).
		\end{aligned}
	\end{equation*}
where $\alpha=max_u|f'(u)|$. This will ensure that
\[
    \frac{\partial f^+(u)}{\partial u}\ge0, \quad\quad
 \frac{\partial f^-(u)}{\partial u}\le0,\quad\quad
 \frac{\partial h^+(u,v)}{\partial v}\ge0, \quad\quad
 \frac{\partial h^-(u,v)}{\partial v}\le0.
\]
And the upwind-biased reconstruction is used for each of those, meaning that a stencil with one more point from the left (resp. right) will be taken for $\{(f^+_j, h^+_j)\}$ (resp. $\{(f^-_j, h^-_j)\}$)
to reconstruct $(\hat{f}^+_{j+1/2}, \hat{h}^+_{j+1/2})$ (resp. $(\hat{f}^-_{j+1/2}, \hat{h}^-_{j+1/2})$),
further giving us that
$$\hat{f}_{i+\frac{1}{2}} =\hat{f}^+_{i+\frac{1}{2}}+ \hat{f}^-_{i+\frac{1}{2}}, \quad 
\hat{h}_{i+\frac{1}{2}} =\hat{h}^+_{i+\frac{1}{2}}+ \hat{h}^-_{i+\frac{1}{2}}.$$ 
In Appendix \ref{sec:app1}, we present the fifth order HWENO construction of $\hat{f}^+_{i+\frac{1}{2}}$ and $\hat{h}^+_{i+\frac{1}{2}}$ given in \cite{fan2023robust}. 
The procedures of $\hat{f}^-_{i+\frac{1}{2}}$ and $\hat{h}^-_{i+\frac{1}{2}}$ are similar by applying a mirror symmetry with respect to $x_{i + 1/2}$.
It is noteworthy that for each grid point $x_i$, we employ a stencil $S=\{(u_{i+j}, v_{i+j}), j=-2, \cdots, +2\}$.

The algorithm of two-dimensional problems or systems is similar to that of the one-dimensional scalar problem described above. Here, we take the two-dimensional system as an example: 
\begin{equation}\label{eq:2DEuler-H}
    \left\{\begin{aligned}
        &\mathbf{U}_t +\mathbf{F}(\mathbf{U})_x +\mathbf{G}(\mathbf{U})_y =\mathbf{0},\\
        &\mathbf{V}_t +\mathbf{F}^x(\mathbf{U})_x +\mathbf{G}^x(\mathbf{U})_y =\mathbf{0},\\
        &\mathbf{W}_t +\mathbf{F}^y(\mathbf{U})_x +\mathbf{G}^y(\mathbf{U})_y=\mathbf{0},
    \end{aligned}\right.
\end{equation}
where
\begin{equation*}
    \begin{aligned}
        &\mathbf{V}=\mathbf{U}_x,  \quad 
        \mathbf{F}^x(\mathbf{U})=\mathbf{F}'(\mathbf{U})\mathbf{V},  \quad 
        \mathbf{G}^x(\mathbf{U})=\mathbf{G}'(\mathbf{U})\mathbf{V}, \\
        &\mathbf{W}=\mathbf{U}_y, \quad 
        \mathbf{F}^y(\mathbf{U})=\mathbf{F}'(\mathbf{U})\mathbf{W}, \quad 
        \mathbf{G}^y(\mathbf{U})=\mathbf{G}'(\mathbf{U})\mathbf{W},
    \end{aligned}
\end{equation*}
and $\mathbf{F}'(\mathbf{U})$ and $\mathbf{G}'(\mathbf{U})$ are the Jacobian matrixes. The semi-discrete finite difference scheme of (\ref{eq:2DEuler-H}) is
\begin{equation}
    \left\{\begin{aligned}
        &\frac{d}{dt}\mathbf{U}_{i,j}=-\frac{1}{\Delta x}(\hat{\mathbf{F}}_{i+\frac{1}{2},j}-\hat{\mathbf{F}}_{i-\frac{1}{2},j})-\frac{1}{\Delta y}(\hat{\mathbf{G}}_{i,j+\frac{1}{2}}-\hat{\mathbf{G}}_{i,j-\frac{1}{2}}), \\
        &\frac{d}{dt}\mathbf{V}_{i,j}=-\frac{1}{\Delta x}(\hat{\mathbf{F}}^x_{i+\frac{1}{2},j}-\hat{\mathbf{F}}^x_{i-\frac{1}{2},j})-\frac{1}{\Delta y}(\hat{\mathbf{G}}^x_{i,j+\frac{1}{2}}-\hat{\mathbf{G}}^x_{i,j-\frac{1}{2}}), \\
        &\frac{d}{dt}\mathbf{W}_{i,j}=-\frac{1}{\Delta x}(\hat{\mathbf{F}}^y_{i+\frac{1}{2},j}-\hat{\mathbf{F}}^y_{i-\frac{1}{2},j})-\frac{1}{\Delta y}(\hat{\mathbf{G}}^y_{i,j+\frac{1}{2}}-\hat{\mathbf{G}}^y_{i,j-\frac{1}{2}}). 
    \end{aligned}\right.
\end{equation}
The numerical fluxes $\hat{\mathbf{F}}_{i\pm\frac{1}{2},j}$ and $\hat{\mathbf{F}}^x_{i\pm\frac{1}{2},j}$ are computed using the one-dimensional HWENO reconstruction procedure applied in the $x$-direction with fixed $j$. 
Similarly, the numerical fluxes $\hat{\mathbf{G}}_{i,j\pm\frac{1}{2}}$ and $\hat{\mathbf{G}}^y_{i,j\pm\frac{1}{2}}$ are obtained through the one-dimensional HWENO reconstruction in the $y$-direction with fixed $i$. 
The remaining numerical fluxes, $\hat{\mathbf{F}}^y_{i\pm\frac{1}{2},j}$ and $\hat{\mathbf{G}}^x_{i,j\pm\frac{1}{2}}$, are derived using the fourth-order linear reconstruction. 
In particular, to mitigate numerical oscillations, the aforementioned reconstruction process is conducted in the characteristic space via characteristic decomposition before being projected back to the physical space.

After the spatial discretization, the third order total variation diminishing (TVD) Runge-Kutta (RK) temporal integration is used for update solution from time $t^n$ to the next level $t^{n+1}=t^n+\Delta t$.
For instance, for the 1D HWENO scheme \eqref{ODSE2}, we have that 
\begin{equation}
    \begin{aligned}
    &\begin{pmatrix}
        u^{(1)}_i\\
        v^{(1)}_i
    \end{pmatrix}
    =
    \begin{pmatrix}
        u^{n}_i\\
        \tilde{v}^n_i
    \end{pmatrix} -\Delta t\begin{pmatrix}
\mathcal{F}_i(u^n,v^n)\\
\mathcal{H}_i(u^n,v^n)
    \end{pmatrix}, \\
    &\begin{pmatrix}
        u^{(2)}_i\\
        v^{(2)}_i
    \end{pmatrix}
    =\frac{3}{4}
    \begin{pmatrix}
        u^{n}_i\\
        \tilde{v}^{n}_i
    \end{pmatrix}+\frac{1}{4} \begin{pmatrix}
        u^{(1)}_i\\
        \tilde{v}^{(1)}_i
    \end{pmatrix} - \frac{1}{4}\Delta t \begin{pmatrix}
        \mathcal{F}_i(u^{(1)},v^{(1)})\\
        \mathcal{H}_i(u^{(1)},v^{(1)})
    \end{pmatrix} , \\
    &\begin{pmatrix}
        u^{n+1}_i\\
        v^{n+1}_i
    \end{pmatrix}
    =
    \frac{1}{3}
    \begin{pmatrix}
        u^{n}_i\\
        \tilde{v}^{n}_i
    \end{pmatrix}+\frac{2}{3}\begin{pmatrix}
        u^{(2)}_i\\
        \tilde{v}^{(2)}_i
    \end{pmatrix} - \frac{2}{3} \Delta t \begin{pmatrix}
        \mathcal{F}_i(u^{(2)},v^{(2)})\\
        \mathcal{H}_i(u^{(2)},v^{(2)})
    \end{pmatrix}.
    \end{aligned}
    \label{RK3}
\end{equation}
We would like to remark that a nonlinear HWENO correction from \( v_i \) to \( \tilde{v}_i \) is necessary. This is because the spatial discretization $\mathcal{H}_i$ is a linear scheme, which can induce numerical oscillations near discontinuities. 
For two-dimensional cases, the HWENO correction is also added at each time step to modify the derivatives in a dimension-by-dimension manner.
The specific procedure for this correction in 1D is detailed in Appendix \ref{sec:app2}.

It is noted in \cite{carpenter1995theoretical} that special care must be taken when we impose the time dependent boundary conditions. Assuming we have the time dependent Dirichlet boundary condition \( g(t) \), the boundary conditions for the intermediate stages of the RK scheme (\ref{RK3}) must be carefully adjusted to maintain the desired third order of accuracy,
\begin{equation}
	\begin{split}
		&u^{n} \sim g(t_n), \\
		&u^{(1)} \sim g(t_n)+\Delta t g'(t_n), \\
		&u^{(2)} \sim g(t_n)+\frac{1}{2}\Delta t g'(t_n)+\frac{1}{4}\Delta t^2 g''(t_n).
	\end{split}
\end{equation}

\section{1D ILW boundary treatment}

At first, we illustrate the idea of the inverse Lax-Wendroff procedure for HWENO reconstruction of 1D scalar conservation laws 
\begin{equation}\label{ODSE}
\begin{cases}
    &u_t+f(u)_x=0, \quad x \in (0,1), \, t>0,\\
    &u(0,t)=g(t), \quad t>0, \\
    &u(x,0)=u_0(t), \quad x\in[0,1].
\end{cases}
\end{equation}
with the assumption $f'(u)>0$. This assumption guarantees the left boundary $x=0$ is an
inflow boundary where a boundary condition should be imposed and the right boundary $x = 1$ is an outflow boundary where no
boundary condition is needed.

Suppose that the grid is divided uniformly as
\begin{equation}
    0\leq C_a\Delta x=x_{1}<x_{2}<...<x_{N}=1-C_b\Delta x\leq1,
\end{equation}
with $C_a, C_b\in[0,1)$. We construct the HWENO scheme \eqref{ODSE2} on each inner grid point, $i=1, \ldots, N$. 
Due to the wide stencil used in HWENO procedure, we need to construct the ghost points values $\{(u_i, v_i)\}$ on $x_{1-i}=x_1-i\Delta x$ and $x_{N+i}=x_N+i\Delta x$, $i=1,2$, with sufficient accuracy. For the outflow boundary condition, we can get the ghost point values via extrapolation directly, or HWENO-type extrapolation to avoid oscillations and maintain accuracy, which will be given later.
Meanwhile, for the inflow boundary condition, we will get them via the ILW boundary treatment in which the given boundary condition is taken into account.

As the ILW procedure outlined in previous work \cite{tan2010inverse,tan2012efficient} operates independently of the interior scheme, we can easily apply it in HWENO scheme directly. The fundamental idea of the ILW procedure is transforming the spatial derivatives at the boundary into time derivatives iteratively utilizing the equations and boundary conditions. Therefore, for the 1D problem (\ref{ODSE}), the ILW procedure will provide us with

\begin{equation} \label{eq:ilw_scalar}
	\begin{split}
		&\partial^{(0)}_xu|_{x=0}=g(t),\\
		&\partial^{(1)}_xu|_{x=0}=-\frac{g'(t)}{f'(g(t))},\\
		&\partial^{(2)}_xu|_{x=0}=\frac{f'(g(t))g''(t)-2f''(g(t))g'(t)^2}{f'(g(t))^3},\\
		&...
	\end{split}
\end{equation}
Then we can construct the 
$d$-th order approximation polynomial through Taylor expansion 
\begin{equation}\label{taylor}
	\tilde{u}(x)=\sum_{k=0}^{d-1}\frac{x^k}{k!}\partial^{(k)}_xu|_{x=0},
\end{equation}
and define 
\begin{equation}
    u_j = \tilde{u}(x_j) = \sum_{k=0}^{d-1} \frac{(x_j)^k}{k!} \partial^{(k)}_x u|_{x=0}, \quad 
    v_j = \tilde{u}'(x_j) = \sum_{k=1}^{d-1} \frac{(x_j)^{k-1}}{(k-1)!} \partial^{(k)}_x u|_{x=0}, \quad 
    j=0,-1. 
\end{equation}

This idea comes from the traditional ILW method for hyperbolic conservation laws proposed by Tan and Shu \cite{tan2010inverse} and will guarantee the global $d$-th order accuracy. 
However, the cost of deriving higher order spatial derivatives becomes prohibitive, especially for multi-dimensional problems and systems. 
To address this, the simplified ILW (SILW) method \cite{tan2012efficient} was proposed when coupling with finite difference schemes, in which the extrapolation is allowed for some higher order derivatives to avoid complex calculations. 
On a stable premise, for a $d$-th order boundary treatment, the ILW process is used to get the spatial derivatives at boundary below order $k_d-1$, while extrapolation is applied for those derivatives of order $k_d$ and above. 
The relation between $d$ and $k_d$ is studied in \cite{li2016stability}.
Recently, Liu et al. \cite{liu2024113259} further reduced $k_d$ by introducing artificial points, thereby improving computational efficiency. 

In fact, we can apply those ideas to construct the ghost point values and derivatives for the HWENO scheme. However, the stability analysis conclusion is suboptimal. For example, for $d=5$, the traditional SILW method requires $k_d=4$ for the fifth-order HWENO scheme, while Liu's modified SILW method requires $k_d=3$. 
In this paper, we improve the Liu's work, allowing the scheme to achieve overall stability with only $k_d=2$.

\subsection{Fifth order SILW boundary treatment.}

\begin{figure}[htb!]
	\centering
	\includegraphics[width=0.7\linewidth]{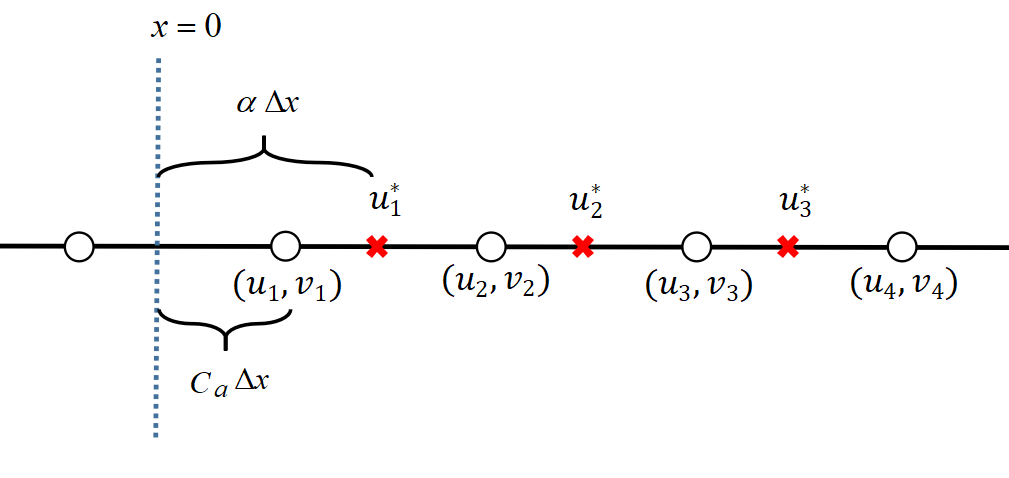}
	\caption{Fifth order SILW boundary treatment diagram.}
	\label{fig:ilw}
\end{figure}

\noindent\textbf{Algorithm 1. Fifth order SILW method for 1D scalar cases.}

\begin{itemize}
\item [Step 1.] Use least square method to construct a quartic polynomial $p(x)$ with inner grid point values $\{(u_i,v_i), 1\leq i\leq k\}$,
\begin{equation*}
    \min_{p(x)\in P^4}\sum_{i=1}^k \left( \|p(x_i)-u_i\|^2 + \Delta x^2\|p'(x_i)-v_i\|^2 \right).
\end{equation*}
To ensure the existence of $p(x)$, we require $k\ge3$.

\item [Step 2.] We use $p(x)$ to construct some special point values 
\[u^*_i=p( i\, \alpha\, \Delta x),\quad 1\leq i\leq 5-k_d.
\]
The parameter $\alpha>0$ will be given later based on the stability analysis.  

\item [Step 3.]
Utilizing the boundary conditions and the ILW procedure to get the spatial derivatives at boundary $\partial_x^{(i)}u|_{x=0}$, $i\leq k_d-1$. 

\item[Step 4.] 
Construct a quartic polynomial $q(x)$ satisfying
\begin{equation*}
    \left\{\begin{aligned}
        &\frac{d^{(i)}}{dx^i} q(0) = \partial_x^{(i)}u|_{x=0},\quad&&0\leq i\leq k_d-1,\\
        &q(i\, \alpha \, \Delta x) = u^*_i,&&1\leq i\leq 5-k_d.
    \end{aligned}\right.
\end{equation*}

\item [Step 5.]
Obtain the values and derivatives at ghost points 
\[
u_{i}=q(x_{i}), \quad v_{i}=q'(x_{i}),\quad
 i=0,-1.\] 
\end{itemize}

Compared with Liu's framework, the main difference is that we introduce an additional parameter \( k \) and incorporate the least squares approach rather than interpolation in Step 1. As the subsequent stability analysis will demonstrate, our method could employ smaller $k_d$, exhibiting superior stability and efficiency.
In the algorithm mentioned above, there are three parameters: $k, k_d$ and $\alpha$. Both $k$ and $k_d$ should be as small as possible for better computational efficiency. However, too small parameters may affect the stability of the algorithm. 
Therefore, we will later determine the selection of $k, k_d$ and $\alpha$ through stability analysis.

\subsection{Fifth order HWENO-type extrapolation}

When the shock wave approaches the boundary, higher order extrapolation may lead to numerical oscillations. To prevent this from happening, we need to design a robust WENO-type extrapolation in Step 4 of Algorithm 1. 
The fifth order WENO-type extrapolation based on $\{ u|_{x=0}, \partial_xu|_{x=0}, u^*_1, u^*_2, u^*_3\}$ in Algorithm 1 will be described in detail.\\

\noindent\textbf{Algorithm 2. Fifth order HWENO-type extrapolation in SILW process.}

\begin{itemize}
\item [Step 1.] Construct five polynomials 
$p_0(x),p_1(x),p_2(x),p_3(x),p_4(x)$ of orders ranging from 0 to 4, satisfying 
\begin{equation*}
    \begin{aligned}
        &p_0(0)=u|_{x=0}, \\
        &p_1(0)=u|_{x=0}, \quad p'_1(0)=\partial_xu|_{x=0},\\
        &p_2(0)=u|_{x=0}, \quad p'_2(0)=\partial_xu|_{x=0}, \quad p_2(i\, \alpha \Delta x)=u^*_i, \;i=1,\\
        &p_3(0)=u|_{x=0}, \quad p'_3(0)=\partial_xu|_{x=0}, \quad p_3(i\, \alpha \Delta x)=u^*_i, \;i=1,2,\\
        &p_4(0)=u|_{x=0}, \quad p'_4(0)=\partial_xu|_{x=0}, \quad p_4(i\, \alpha \Delta x)=u^*_i, \;i=1,2,3.
    \end{aligned}
\end{equation*}
And define their linear weights $d_r$ to be
\begin{equation*}
    d_0=\Delta x^4,\;d_1=\Delta x^3,\;d_2=\Delta x^2,\;d_3=\Delta x,\;d_4=1-\sum_{i=0}^3 d_i.
\end{equation*}

\item[Step 2.] Define the smoothness indicators 
\begin{equation}
\beta_i=\sum_{j=1}^{i}\int^{\Delta x/2}_{-\Delta x/2}\Delta x^{2j-1}(\frac{d^{j}p_i(x)}{dx^{j}})^2dx,\quad i=1,2,3,4, 
\end{equation}
which is used to measure the smoothness of each polynomial. To ensure $\beta_i$ at the same magnitude in the smooth case, we take $\beta_0=\Delta x^2$.

\item[Step 3.] 
Next, express the nonlinear weight $\omega_i$ as 
\begin{equation*}
    \omega_i=\frac{\gamma_i}{\sum_{j=0}^{4}\gamma_j},\quad\gamma_i=\frac{d_i}{(\epsilon+\beta_i)^2}.
\end{equation*}
Here, $\epsilon$ is used to avoid the denominator being zero, and it is typically set to $\epsilon=10^{-6}$.

\item[Step 4.]
Finally, we can ultimately obtain the extrapolation polynomial $p(x)$
\begin{equation*}
    p(x)=\sum_{i=0}^{4} \omega_i p_i(x).
\end{equation*}

\end{itemize}

With the help of smoothness indicator and nonlinear weights, the extrapolation can achieve fifth order accuracy when the solution is smooth near boundary, and degenerate to lower order scheme to avoid spurious oscillation if there is a discontinuity close to the boundary.

For the least squares method, we can use a similar HWENO-type extrapolation. A constant polynomial \( p_0(x) \) is constructed from \( \{u_0\} \). A quadratic polynomial \( p_1(x) \) is obtained by least squares fitting \( \{u_1, v_1, u_2, v_2\} \), and a quartic polynomial \( p_2(x) \) is obtained by least squares fitting \( \{u_1, v_1, u_2, v_2, u_3, v_3\} \). The linear weights are \( (d_0, d_1, d_2) = (\Delta x^4, \Delta x^2, 1 - \Delta x^2 - \Delta x^4) \). The nonlinear weights \( \omega_0, \omega_1, \omega_2 \) are computed as above, giving the extrapolation polynomial \( p(x) = \omega_0 p_0(x) + \omega_1 p_1(x) + \omega_2 p_2(x) \).

On the other hand, the HWENO-type extrapolation can be applied for the outflow boundary $x=b$ as well. We construct a linear polynomial \( p_0(x) \) using the interior points \( \{u_N, v_N\} \), a cubic polynomial \( p_1(x) \) using \( \{u_N, v_N, u_{N-1}, v_{N-1}\} \), and a fifth-degree polynomial \( p_2(x) \) using \( \{u_N, v_N, u_{N-1}, v_{N-1}, u_{N-2}, v_{N-2}\} \). The linear weights are assigned as \( (d_0, d_1, d_2) = (\Delta x^4, \Delta x^2, 1 - \Delta x^2 - \Delta x^4) \). Finally, the nonlinear weights \( \omega_0, \omega_1, \omega_2 \) are constructed using the same method as described above, yielding the extrapolation polynomial \( p(x)= \omega_0 p_0(x) + \omega_1 p_1(x) + \omega_2 p_2(x) \). Thus, we can obtain $u_{N+i}=p(x_{N+i})$ and $v_{N+i}=p'(x_{N+i})$, $i=1,2$.

\subsection{Stability analysis}

Based on the advection equation $u_t+cu_x=0$, we analyze the stability of the fully discrete scheme, employing the fifth order linear upwind scheme, third order TVD Runge-Kutta time evolution and the proposed boundary treatment. 
Assume $c>0$, thus the upwind linear scheme can be expressed as
\begin{equation}
\begin{aligned}
    \frac{du_i}{dt}=&-\frac{c}{\Delta x} \left( \frac{23}{120}u_{i-2}-\frac{33}{40}u_{i-1}+\frac{3}{40}u_i+\frac{67}{120}u_{i+1}+\frac{3\Delta x}{40}v_{i-2}-\frac{3\Delta x}{40}v_{i-1}+\frac{7\Delta x}{40}v_i-\frac{7\Delta x}{40}v_{i+1} \right), \\
    \frac{dv_i}{dt}=& -\frac{c}{\Delta x}\left( -\frac{3}{8\Delta x}u_{i-2}+\frac{19}{8\Delta x}u_{i-1}-\frac{29}{8\Delta x}u_i+\frac{13}{8\Delta x}u_{i+1}-\frac{1}{8}v_{i-2}+\frac{1}{8}v_{i-1}+\frac{3}{8}v_i-\frac{3}{8}v_{i+1} \right).
\end{aligned}
\label{space}
\end{equation}

Let us start with the periodic boundary problems. Substituting $u(x,t)=\frac{1}{\sqrt{2\pi}} e^{i\omega x}\hat{u}(\omega,t)$ and $v(x,t)=\frac{1}{\sqrt{2\pi}} e^{i\omega x}\hat{v}(\omega,t)$ into the above equations gives us a matrix-vector form
\begin{equation}
	\frac{d\hat{\mathbf{u}}}{dt}=\frac{c}{\Delta x}\mathcal{Q}_1\hat{\mathbf{u}}.
\end{equation}
Here $\hat{\mathbf{u}}=(\hat{u}, \hat{v})^T$. The matrix $\mathcal{Q}_1$ can be expanded as 
\begin{equation*}
	\mathcal{Q}_1=\left(
	\begin{matrix}
		q_1 & q_2\\
		q_3 & q_4
	\end{matrix}\right),
\end{equation*}
with $\xi = \omega\Delta x$ and 
\begin{equation*}
    \left\{\begin{aligned}
        &q_1=\left(-\frac{23}{120}\cos2\xi+\frac{4}{15}\cos\xi-\frac{3}{40}\right)
        +\left(\frac{23}{120}\sin2\xi-\frac{83}{60}\sin\xi \right)i, \\
        &q_2=\Delta x \left[ \left(-\frac{3}{40}\cos2\xi+\frac{1}{4}\cos\xi-\frac{7}{40}\right)
        +\left(\frac{3}{40}\sin2\xi+\frac{1}{10}\sin\xi\right)i \right], \\
        &q_3=\frac{1}{\Delta x}\left[ \left(\frac{3}{8}\cos2\xi-4\cos\xi+\frac{29}{8}\right)
        +\left(-\frac{3}{8}\sin2\xi+\frac{3}{4}\sin\xi\right)i \right], \\
        &q_4=\left(\frac{1}{8}\cos2\xi+\frac{1}{4}\cos\xi-\frac{3}{8}\right)
        +\left(-\frac{1}{8}\sin2\xi+\frac{1}{2}\sin\xi\right)i .
    \end{aligned}\right.
\end{equation*}

It is noted that the HWENO scheme uses  the HWENO limiter to modify the first-order derivative $v_i$ in each RK stage to suppress oscillations. If we consider the linear modification, then $\tilde{v}_i$ can be expressed as
\begin{equation}
    \tilde{v}_i=\frac{3}{4\Delta x}(u_{i+1}-u_{i-1})-\frac{1}{4}(v_{i+1}+v_{i-1}).
    \label{modify}
\end{equation}
The corresponding modification matrix form is
\begin{equation*}
    \begin{pmatrix}
        \hat{u} \\ \hat{\tilde{v}}
    \end{pmatrix} = 
    \mathcal{Q}_2 \begin{pmatrix}
        \hat{u} \\ \hat{v}
    \end{pmatrix} ,\qquad
    \mathcal{Q}_2=\left(\begin{matrix}
        1 & 0\\
        \frac{3}{2\Delta x}\sin\xi i& -\frac{1}{2}\cos\xi
    \end{matrix}\right).
\end{equation*}

\begin{figure}[htb!]
    \centering
    \subfigure[$\lambda_{cfl}=1.07$]{\includegraphics[width=1.0in,angle=0,scale=3.0]{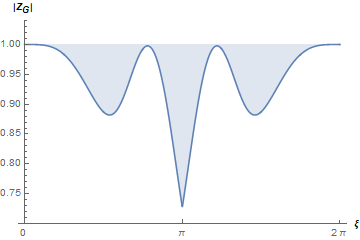}}
    \subfigure[$\lambda_{cfl}=1.08$]{\includegraphics[width=1.0in,angle=0,scale=3.0]{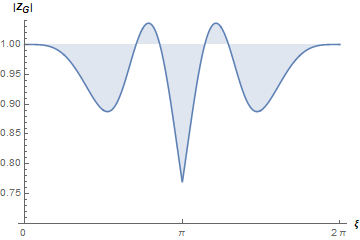}}
    \caption{The absolute value of the maximum eigenvalue $|z_G|$ with respect to $\xi$. }
    \label{fig:stable_perodic}
\end{figure}

Combined with the third order TVD RK time evolution \eqref{RK3} and denote $\lambda_{cfl}=\frac{c\Delta t}{\Delta x}$,
we have that 
$$\hat{\mathbf{u}}^{n+1} = \mathcal{G} \hat{\mathbf{u}}^n, $$
and the matrix $\mathcal{G}$ can be written as follows
\begin{equation}
    \mathcal{G}=\frac{1}{3}\mathcal{Q}_2+\frac{2}{3}(\mathcal{Q}_2+\lambda_{cfl}\mathcal{Q}_1)(\frac{3}{4}\mathcal{Q}_2+\frac{1}{4}(\mathcal{Q}_2+\lambda_{cfl}\mathcal{Q}_1)(\mathcal{Q}_2+\lambda_{cfl}\mathcal{Q}_1)).
\end{equation}
In order to ensure the stability of the fully discrete scheme, the two eigenvalues of the matrix $\mathcal{G}$ ($z_1$ and $z_2$) have to fall inside the unit circle, i.e., $|z_1|\leq1,|z_2|\leq1$ for all $\xi$. 
It is worth noting that the $\Delta x$ and $\frac{1}{\Delta x}$ will cancel each other out when calculating the eigenvalues, resulting that $z_1$ and $z_2$ are independent of $\Delta x$. It is easy to prove that the stability condition is  $\lambda_{cfl}\leq1.07$. We plot the absolute value of the maximum eigenvalue $|z_G|=\max(|z_1|, |z_2|)$ in Figure \ref{fig:stable_perodic} to verify the CFL condition. 

Next, we will concentrate on our methods for solving the linear advection equation involving homogeneous Dirichlet conditions $g(t) = 0$ at $x=0$. Our starting point is to ensure that the boundary treatment will not affect the stability of inner scheme, that is the scheme can maintain stable under $\lambda_{cfl}\leq 1.07$ for any $C_a\in[0,1)$.

Initially, we express the semi-discrete scheme (\ref{space}) with the proposed numerical boundary treatments in matrix-vector form
\begin{equation}
    \frac{d}{dt}\mathbf{U}=\frac{c}{\Delta x}\mathbf{Q}_1\mathbf{U},
\end{equation}
where, $\mathbf{U}=(u_1, \cdots, u_N,v_1, \cdots, v_N)^T$, and $\mathbf{Q}_1$ is the coefficient matrix corresponding to the spatial discretization and boundary treatment. At the same time, we can obtain the modification matrix of the corresponding spatial discretization (\ref{modify})
\begin{equation}
  \tilde{\mathbf{U}}=\mathbf{Q}_2\mathbf{U}.
\end{equation}
Combined with the time evolution (\ref{RK3}), the fully discrete can be presented as 
\begin{equation}
\mathbf{U}^{n+1} =  \mathbf{G} \mathbf{U}^n,
\end{equation}
with the matrix 
\begin{equation}
    \mathbf{G}=\frac{1}{3}\mathbf{Q}_2+\frac{2}{3}(\mathbf{Q}_2+\lambda_{cfl}\mathbf{Q}_1)(\frac{3}{4}\mathbf{Q}_2+\frac{1}{4}(\mathbf{Q}_2+\lambda_{cfl}\mathbf{Q}_1)(\mathbf{Q}_2+\lambda_{cfl}\mathbf{Q}_1)).
\end{equation}

Relevant studies \cite{vilar2015development,li2016stability} indicate that to ensure the fully discrete scheme is stable, we only need to focus on the fixed eigenvalues of the matrix $\mathbf{G}$ that do not change with respect to the number of grids $N$, i.e., $|z_\mathbf{G}|\leq 1$ holds for $\lambda_{cfl} \leq 1.07$ and all $C_a\in[0,1)$. 
Furthermore, it would be even better if the algorithm is stable with $k_d$ and $k$ as less as possible, as this not only ensures a reduced computational cost, but also matches the compactness of HWENO scheme.

Table \ref{tab:stablity} shows the range of parameter $\alpha$ in our algorithm under the stability conditions for different $k$ and $k_d$. 
For example, when \( k = 3 \) and \( k_d = 1 \), no parameter \( \alpha \) exists to ensures the the algorithm stable for all \( C_a \in [0, 1) \). 
When \( k = 3 \) and \( k_d = 2 \), the algorithm is stable if \( \alpha \) is chosen within the interval \([0.93, 1.09]\). Therefore, in the one-dimensional cases, we select \( k = 3 \) and \( k_d = 2 \) to implement our algorithm.

\begin{table}[htb!]
	\centering
	\caption{Linear stability analysis results: reasonable range of parameter $\alpha$.}
	\begin{tabular}{|c|c|c|}\hline
		\diagbox{$k_d$}{$k$}& 3   &4   \\ \hline
		$k_d$=1 & - & - \\ \hline
		$k_d$=2 & [0.93,1.09] & [0.84,1.74] \\ \hline
		$k_d$=3 & [0.77,10+] & [0.78,10+]  \\ \hline
	\end{tabular}
\label{tab:stablity}
\end{table}

\section{2D ILW boundary treatment}

Consider the two-dimensional compressible Euler equations along with suitable boundary and initial conditions. The governing equation is given as
\begin{equation}\label{2DEuler}
\mathbf{U}_t+\mathbf{F}(\mathbf{U})_x+\mathbf{G}(\mathbf{U})_y=\mathbf{0}, \quad (x,y)\in \Omega, 
\end{equation}
with
\begin{equation}
	\mathbf{U} = \begin{pmatrix}
		\rho\\
		\rho u\\
		\rho v\\
		E
	\end{pmatrix}, \quad 
	\mathbf{F} = \begin{pmatrix}
		\rho u\\
		\rho u^2+p\\
		\rho uv\\
		u(E+p)
	\end{pmatrix}, \quad
	\mathbf{G} = \begin{pmatrix}
		\rho v\\
		\rho uv\\
		\rho v^2+p\\
		v(E+p)
	\end{pmatrix}, 
\end{equation}
and 
\begin{equation*}
	E=\frac{p}{\gamma-1}+\frac{1}{2}\rho(u^2+v^2).
\end{equation*}
Here, $ \rho $ is the density, $\mathbf{u}=(u,v)^T$ is the velocity, $ p $ is the pressure, $ E $ is the total energy, and $\gamma=1.4$ for ideal polytropic gas.

The domain $\Omega$ can be arbitrary complex geometry. We assume the domain is covered by a uniform non-body-fitting Cartesian mesh, with 
$$\Delta x=x_{i+\frac{1}{2}}-x_{i-\frac{1}{2}},\quad
\Delta y=y_{j+\frac{1}{2}}-y_{j-\frac{1}{2}}, \quad 
h=\sqrt{\Delta x^2+\Delta y^2}.$$
Suppose $\Delta x$ and $\Delta y$ are of the same magnitude.

Since the HWENO method requires information from neighboring points, we need to construct information for ghost points that fall outside the domain. 
Assuming that $P_{i,j}$ is a ghost point. For the fifth order HWENO scheme, we need to construct the point value $\mathbf{U}_{i,j}$ with fifth order accurate and the first derivatives $(\mathbf{U}_x)_{i,j}, (\mathbf{U}_y)_{i,j}$ with fourth order accurate via SILW method.

\begin{figure}[htb!]
	\centering
	\includegraphics[width=0.5\linewidth]{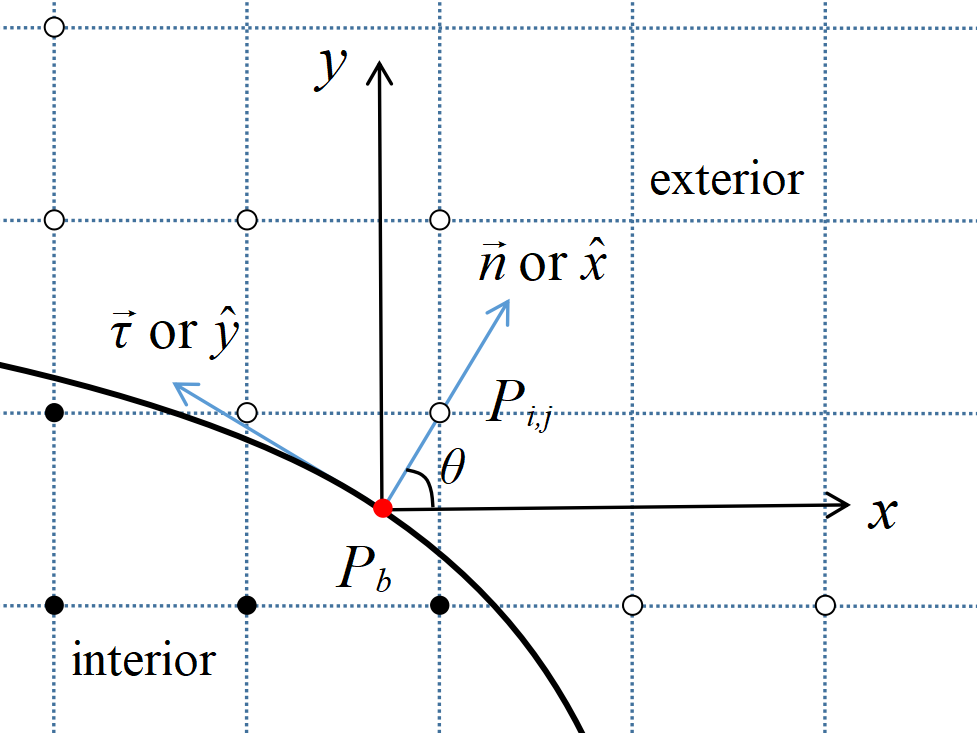}
	\caption{Two dimensional ILW method diagram}
	\label{fig:2dilw}
\end{figure}

Similar to \cite{tan2010inverse, tan2012efficient, liu2024113259}, at first, we need to perform a local coordinate transformation to convert the 2D extrapolation into 1D cases. Find the corresponding foot point $P_b\in\partial\Omega$ whose outward normal vector passes through $P_{i,j}$, and denote the unit outward normal vector as $\mathbf{n}=(\cos\theta, \sin\theta)$, see Figure \ref{fig:2dilw}. 
At point $P_b$, we set up a local coordinate system such that the $\hat{x}$-axis aligns with the direction of $\mathbf{n}$ and the $\hat{y}$-axis is in the tangential direction of $\partial\Omega$, i.e.,
\begin{equation}\label{eq:rotation}
	\begin{pmatrix}
		\hat{x}\\ \hat{y}
	\end{pmatrix}
	=
	\begin{pmatrix}
		\cos\theta&\sin\theta\\
		-\sin\theta&\cos\theta
	\end{pmatrix}
	\begin{pmatrix}
		x\\ y
	\end{pmatrix}.
\end{equation}
In the new coordinate system, the equation (\ref{2DEuler}) can be expressed as
\begin{equation}\label{2DEuler_2}
    \hat{\mathbf{U}}_t+\mathbf{F}(\hat{\mathbf{U}})_{\hat{x}}+\mathbf{G}(\hat{\mathbf{U}})_{\hat{y}}=\mathbf{0}, 
\end{equation}
with
\begin{equation*}
	\hat{\mathbf{U}} = \begin{pmatrix}
		\rho\\
		\rho \hat{u}\\
		\rho \hat{v}\\
		E
	\end{pmatrix} = 
    \begin{pmatrix}\hat{U}_1\\ \hat{U}_2 \\ \hat{U}_3 \\ \hat{U}_4\end{pmatrix},\qquad 
	\begin{pmatrix}
		\hat{u}\\
		\hat{v}
	\end{pmatrix}=
	\begin{pmatrix}
		\cos\theta&\sin\theta\\
		-\sin\theta&\cos\theta
	\end{pmatrix}
	\begin{pmatrix}
		u\\
		v
	\end{pmatrix}.
\end{equation*}
Let $\mathbf{A}(\mathbf{\hat{U}})$ be the Jacobian matrix of $\mathbf{F}(\hat{\mathbf{U}})_{\hat{x}}$,
\begin{equation*}
	\mathbf{A}(\mathbf{\hat{U}}) =\mathbf{F}'(\mathbf{\hat{U}}) 
        =\begin{pmatrix}
		\mathbf{a}_1(\mathbf{\hat{U}})\\
		\mathbf{a}_2(\mathbf{\hat{U}})\\
		\mathbf{a}_3(\mathbf{\hat{U}})\\
		\mathbf{a}_4(\mathbf{\hat{U}})
	\end{pmatrix}.
\end{equation*}
Since (\ref{2DEuler_2}) is hyperbolic, the matrix $\mathbf{A}(\mathbf{\hat{U}})$ is diagonalized
$\mathbf{A}(\mathbf{\hat{U}}) =\mathbf{L}^{-1} \mathbf{\Lambda} \mathbf{L}$
with
\begin{equation*}
    \mathbf{\Lambda}(\mathbf{\hat{U}})=diag\{\hat{u}-c,\hat{u},\hat{u},\hat{u}+c\},\quad
    c=\sqrt{\gamma p/\rho}, \quad
    \mathbf{L}(\mathbf{\hat{U}})
        =\begin{pmatrix}
		\mathbf{l}_1(\mathbf{\hat{U}})\\
		\mathbf{l}_2(\mathbf{\hat{U}})\\
		\mathbf{l}_3(\mathbf{\hat{U}})\\
		\mathbf{l}_4(\mathbf{\hat{U}})
	\end{pmatrix}.
\end{equation*}

The number of boundary conditions required at point $P_b$ depends on the sign of eigenvalues $\{\hat{u}-c,\hat{u},\hat{u},\hat{u}+c\}$.
\begin{itemize}
	\item[Case 1.] $\hat{u}-c\geq 0$: no boundary condition is needed.
	\item[Case 2.] $\hat{u}-c<0 \leq \hat{u}$: only one boundary condition is needed.
	\item[Case 3.] $\hat{u}<0 \leq\hat{u}+c$: three boundary conditions are needed.
	\item[Case 4.] $\hat{u}+c<0$: four boundary conditions are needed.
\end{itemize}

In the following, we will take case 3 as an example to describe our SILW procedure, assuming that the given boundary conditions are Dirichlet boundary conditions at $P_b$,
\begin{equation*}
    \hat{U}_1=g_1(t),\quad
    \hat{U}_2=g_2(t),\quad
    \hat{U}_3=g_3(t).
\end{equation*}

For the 2D problem, we will apply the 1D algorithm along $\hat{x}$-direction with the help of the local coordinate transformation.
We should point out that, suppose in 1D we take $\alpha\in[\alpha_1, \alpha_2]$ for stability, it should be 
\(\alpha\in[\alpha_1, \alpha_2 \frac{\min(\Delta x, \Delta y)}{h}]\) in 2D. 
In order to ensure the interval is not null, for 2D cases
we choose \( k = 4 \) and \( k_d = 2 \),  
which provides a larger stability region. Thus, we need to construct values on three artificial points, and
the following values via ILW procedure
\begin{equation}
    \hat{\mathbf{U}}^* = \hat{\mathbf{U}}|_{P_b} +\mathcal{O}(h^5), \quad
    \hat{\mathbf{U}}_{\hat{x}}^* = \hat{\mathbf{U}}_{\hat{x}}|_{P_b} +\mathcal{O}(h^4), \quad
    \hat{\mathbf{U}}_{\hat{y}}^* = \hat{\mathbf{U}}_{\hat{y}}|_{P_b} +\mathcal{O}(h^4).
\end{equation}
In particular, with the local characteristic decomposition, we will use the ILW procedure for the ingoing characteristic variables and extrapolation for the outgoing characteristic variables.

In the local characteristic decomposition, we need to give the approximation of $\hat{\mathbf{U}}$ at $P_b$ at first.
Define a quartic polynomial vector $\mathbf{P}(x,y)$ as the approximating polynomial obtained through least squares based on internal point values and derivatives with sufficient grid points, and we can obtain $\hat{\mathbf{P}}(\hat{x},\hat{y})$ through the rotation transformation (\ref{eq:rotation}). 
Then, we can define $\mathbf{\hat{U}}_{ext}= \hat{\mathbf{P}}|_{P_b}$ and perform a local characteristic decomposition $\mathbf{V}=\mathbf{L}(\mathbf{\hat{U}}_{ext}) \mathbf{\hat{U}} =(V_1,V_2, V_3, V_4)^T$ for the points near the point $P_b$.
Note that $V_4$ is the outgoing characteristic variable. Thus, we can get its values and derivatives at $P_b$ by extrapolation (or HWENO extrapolation) directly.

Combining the given boundary conditions and extrapolation on $V_4$, we can obtain the following system of equations,
\begin{equation}\label{eq:ILW_2D_0}
        \left\{
	\begin{aligned}
		&\hat{U}_1^{*} = g_1(t),\\
		&\hat{U}_2^{*} = g_2(t),\\
        &\hat{U}_3^{*} = g_3(t),\\
		&\mathbf{l}_4(\mathbf{\hat{U}}_{ext}) \cdot \mathbf{\hat{U}}^{*} = V_4^{*},
	\end{aligned}\right.
\end{equation}
where $V_4^{*}=\mathbf{l}_4 (\mathbf{\hat{U}}_{ext})\cdot \mathbf{\hat{U}}_{ext}$. Solving the above system gives us the point value $\mathbf{\hat{U}}^*$.

Next we try to find the spatial derivatives $\mathbf{\hat{U}}^*_{\hat{x}}$ by the ILW procedure. We can place the term $\mathbf{G}(\hat{\mathbf{U}})_{\hat{y}}$ in (\ref{2DEuler_2}) as a source term on the right hand side, 
\begin{equation}\label{2DNCF}
	\hat{\mathbf{U}}_t + \mathbf{A}(\hat{\mathbf{U}}) \hat{\mathbf{U}}_{\hat{x }}=\mathbf{Resy}, \quad
    \mathbf{Resy}=-\frac{\partial \mathbf{G}(\mathbf{\hat{U}})}{\partial \hat{y}}=\begin{pmatrix}
		Resy_1\\
		Resy_2\\
		Resy_3\\
		Resy_4
	\end{pmatrix}.
\end{equation}
Plugging the boundary conditions in the first three equations of \eqref{2DNCF} and combining with the extrapolated values $(V_4^{*})_{\hat{x}}$ gives us that
\begin{equation}\label{eq:ILW_2D_1}
        \left\{
	\begin{aligned}
		&\mathbf{a}_1(\mathbf{\hat{U}}^*)\cdot \mathbf{\hat{U}}_{\hat{x}}^*=-g'_1(t)+Resy_1, \\
		&\mathbf{a}_2(\mathbf{\hat{U}}^*)\cdot \mathbf{\hat{U}}_{\hat{x}}^*=-g'_2(t)+Resy_2, \\
        &\mathbf{a}_3(\mathbf{\hat{U}}^*)\cdot \mathbf{\hat{U}}_{\hat{x}}^*=-g'_3(t)+Resy_3, \\
		&\mathbf{l}_4(\mathbf{\hat{U}}_{ext}) \cdot \mathbf{\hat{U}}_{\hat{x}}^{*} = (V_4)^{*}_{\hat{x}}.
	\end{aligned}
 \right.
\end{equation}
Here, all values in $Resy_1,Resy_2,Resy_3$ and the derivatives of the outgoing variable $(V_4)^{*}_{\hat{x}}$ are obtained by extrapolation polynomial $\hat{\mathbf{P}}$. 
Hence, we can get $\mathbf{\hat{U}}^{*}_{\hat{x}}$ by solving the above system.

Similarly, equation (\ref{2DEuler_2}) can also be rewritten as
\begin{equation}
	\hat{\mathbf{U}}_t + \mathbf{B}(\hat{\mathbf{U}}) \hat{\mathbf{U}}_{\hat{y}} =\mathbf{Resx}, 
\end{equation}
with
\begin{equation*}
    \mathbf{B}(\hat{\mathbf{U}}) = \mathbf{G}'(\hat{\mathbf{U}}) =\begin{pmatrix}
		\mathbf{b}_1(\hat{\mathbf{U}})\\
		\mathbf{b}_2(\hat{\mathbf{U}})\\
		\mathbf{b}_3(\hat{\mathbf{U}})\\
		\mathbf{b}_4(\hat{\mathbf{U}})
	\end{pmatrix}, \quad 
    \mathbf{Resx}=-\frac{\partial \mathbf{F}(\mathbf{\hat{U}})}{\partial \hat{x}}=\begin{pmatrix}
		Resx_1\\
		Resx_2\\
		Resx_3\\
		Resx_4
	\end{pmatrix}.
\end{equation*}
Then, $\mathbf{\hat{U}}^{*}_{\hat{y}}$ can be obtained by solving the linear system
\begin{equation}\label{eq:ILW_2D_2}
        \left\{
	\begin{aligned}
		&\mathbf{b}_1(\mathbf{\hat{U}}^*)\cdot \mathbf{\hat{U}}_{\hat{y}}^*=-g'_1(t) +Resx_1,\\
        &\mathbf{b}_2(\mathbf{\hat{U}}^*)\cdot \mathbf{\hat{U}}_{\hat{y}}^*=-g'_2(t) +Resx_2,\\
        &\mathbf{b}_3(\mathbf{\hat{U}}^*)\cdot \mathbf{\hat{U}}_{\hat{y}}^*=-g'_3(t) +Resx_3,\\
		&\mathbf{l}_4(\mathbf{\hat{U}}_{ext}) \cdot \mathbf{\hat{U}}_{\hat{y}}^{*} = (V_4)^{*}_{\hat{y}}.
	\end{aligned}
 \right.
\end{equation}
All values in $Resx_1,Resx_2,Resx_3$ and the derivatives $(V_4)^{*}_{\hat{y}}$ are obtained by extrapolation as well. 

Once we have $\mathbf{\hat{U}}^{*}$, $\mathbf{\hat{U}}^{*}_{\hat{x}}$ and $\mathbf{\hat{U}}^{*}_{\hat{y}}$, we can construct the extrapolation along $\hat{x}$-direction. 
In particular, we will construct a quartic function $\mathbf{Q}_0(s) = \hat{\mathbf{U}}(P_b +s \mathbf{n})+\mathcal{O}(h^5)$, and a cubic function $\mathbf{Q}_1(s) = \hat{\mathbf{U}}_{\hat{y}}(P_b +s \mathbf{n})+\mathcal{O}(h^4)$, which will be employed to
evaluate $\{\mathbf{\hat{U}}^{*}, \mathbf{\hat{U}}^{*}_{\hat{x}} \}$ and $\mathbf{\hat{U}}^{*}_{\hat{y}}$ at $P_{i,j}$, respectively.\\

\noindent
\textbf{Algorithm 3. Fifth order SILW method for 2D Euler equations}

\begin{itemize}
\item [Step 1.]For the ghost point $P_{ij}$ outside the computational domain, we drop a perpendicular to find the foot of the perpendicular $P_b$ on the boundary. We then take the internal point values and derivative values near $P_b$ to perform a least squares fitting to construct $\mathbf{P}(x,y)$. Then, we apply a local coordinate transformation resulting in $\hat{\mathbf{P}}(\hat{x},\hat{y})$.

\item [Step 2.]
Construct point value $\hat{\mathbf{U}}^{*}$ and spatial derivatives $(\hat{\mathbf{U}}^{*})_{\hat{x}},(\hat{\mathbf{U}}^{*})_{\hat{y}}$ at $P_b$ via the ILW procedure by solving \eqref{eq:ILW_2D_0}, \eqref{eq:ILW_2D_1} and \eqref{eq:ILW_2D_2}.

\item [Step 3.] Construct informations on three artificial inner points
\begin{equation*}
    \begin{aligned}
        \mathbf{Z}_{k}^{(0)}=\hat{\mathbf{P}}(P_b-kh\mathbf{n}),\quad 
        \mathbf{Z}_{k}^{(1)}=\hat{\mathbf{P}}_{\hat{y}}(P_b-kh\mathbf{n}),\quad k=1,2,3.
    \end{aligned}
\end{equation*}

\item [Step 4.]Construct a quartic polynomial $\mathbf{Q}_0(s)$ and a cubic polynomial $\mathbf{Q}_1(s)$ satisfying
\begin{align*}
    &\left\{
        \begin{aligned}
        &\mathbf{Q}_0(0)=\hat{\mathbf{U}}^{*},\quad
        \mathbf{Q}'_0(0)=(\hat{\mathbf{U}}^{*})_{\hat{x}},\\
        &\mathbf{Q}_0(-kh)=\mathbf{Z}_{k}^{(0)},\quad k=1,2,3,
        \end{aligned}
    \right.\\
    &\left\{
        \begin{aligned}
        &\mathbf{Q}_1(0)=(\hat{\mathbf{U}}^{*})_{\hat{y}},\\
        &\mathbf{Q}_1(-kh)=\mathbf{Z}_{k}^{(1)},\quad k=1,2,3.
        \end{aligned}
    \right.
\end{align*}
\item [Step 5.] Take $\bar{s}=|P_{i,j}-P_b|$. Thus, we can obtain the point values and derivative values at $P_{i,j}$ using  $\mathbf{Q}_0$ and $\mathbf{Q}_1$,
$$\hat{\mathbf{U}}_{i,j} = \mathbf{Q}_0(\bar{s}), \quad 
(\hat{\mathbf{U}}_{\hat{x}})_{i,j} = \mathbf{Q}'_0(\bar{s}),\quad 
(\hat{\mathbf{U}}_{\hat{y}})_{i,j} = \mathbf{Q}_1(\bar{s}).$$
Then, we can obtain $\mathbf{U}_{i,j}$, $(\mathbf{U}_x)_{i,j}$ and $(\mathbf{U}_y)_{i,j}$ by applying a rotation transformation back to $(x, y)$. 
\end{itemize}

\section{Numerical tests}

In this section, we will validate the effectiveness and stability of our SILW method through several numerical experiments. We will employ the fifth-order HWENO-R scheme for spatial discretization, coupled with the SILW method for boundary treatment. Specifically, for one-dimensional tests, we adopt \( (k_d, k) = (2, 3) \), while for two-dimensional tests, we use \( (k_d, k) = (2, 4) \). For all cases, we set the parameter \( \alpha = 1 \). Time discretization is performed using a third-order RK scheme (\ref{RK3}). The time step is chosen as:
\[
\Delta t = 0.6 \frac{\Delta x^{k}}{a_x}
\]
for one-dimensional problems, and
\[
\Delta t = \frac{0.6}{a_x / \Delta x^{k} + a_y / \Delta y^{k}}
\]
for two-dimensional problems. Here, \( a_x = \max |\lambda (\mathbf F'(\mathbf U))| \), \( a_y = \max |\lambda (\mathbf G'(\mathbf u))| \), where \( \lambda \) represents the eigenvalues of the Jacobian matrix. 
To ensure fifth-order accuracy, $k=5/3$ for the accuracy tests. Otherwise, we just take $k=1$.

\subsection{One dimensional cases}

\vspace{0.3cm}
\noindent
\textbf{Example 1.}
At first, we consider the 1D Burgers' equation
 \begin{equation*}
    \begin{cases}
        u_t+(\frac{u^2}{2})_x=0,\quad &x\in(0,2), \;\; t>0,\\
        u(x,0)=1+\sin(\pi x), \quad &x\in[0,2],\\
        u(0,t)=g(t),&t>0.
    \end{cases}
\end{equation*}
The left boundary $x=0$ is an inflow boundary, where the boundary condition $g(t)$ is precisely solved using the Newton iteration method for the problem with periodic boundaries, and the time derivative at the boundary is obtained through interpolation. The right boundary $x=2$ is an outflow boundary, where no boundary condition is needed. 
We divide the domain with uniform mesh size $\Delta x=2/N$ and set $C_b = 1-C_a$.

\begin{figure}[htb!]
    \centering
    \subfigure[$T=1.5/\pi$]
    {\includegraphics[width=1.0in,angle=0,scale=3.0]{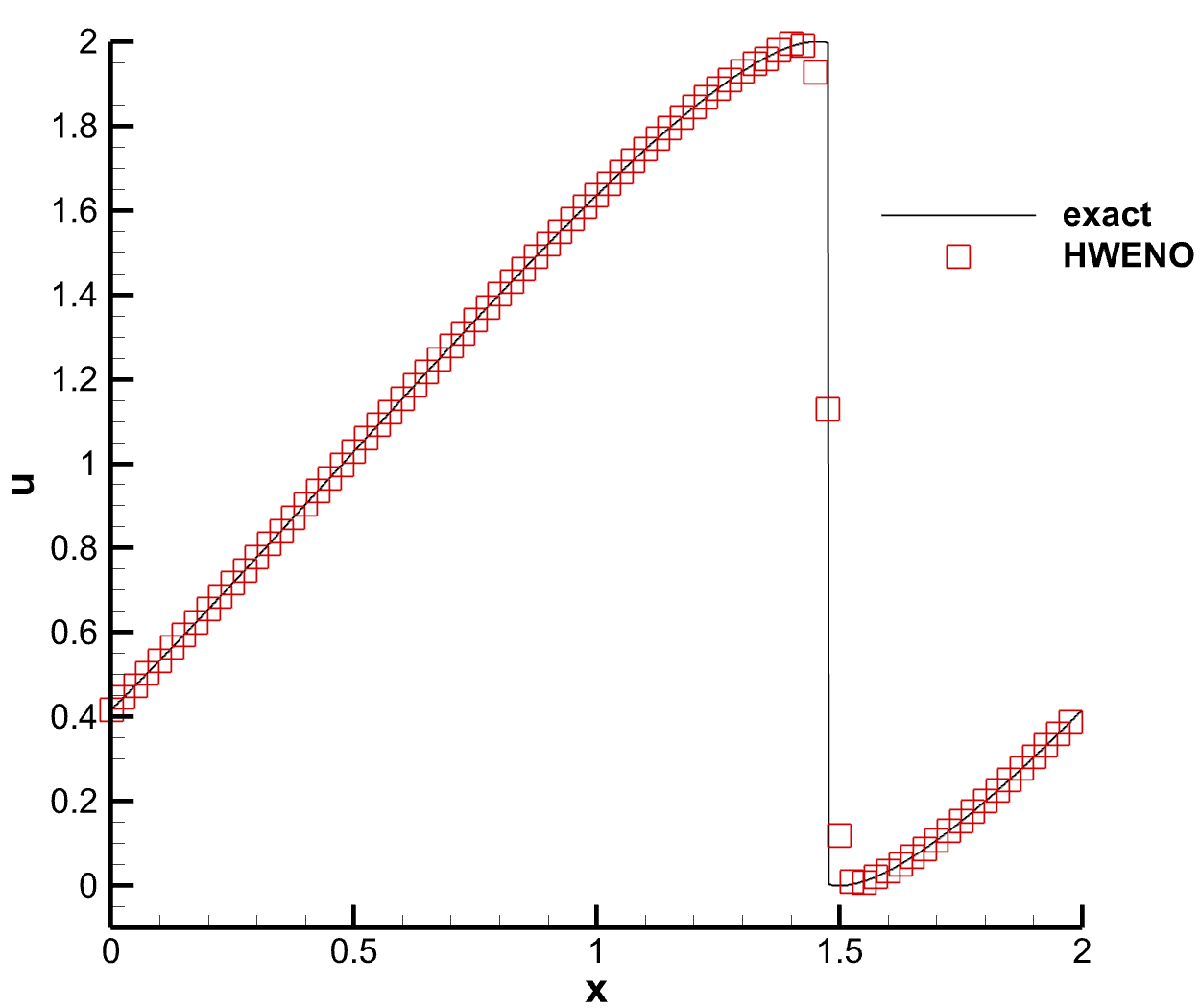}}
    \subfigure[$T=5/\pi$]
    {\includegraphics[width=1.0in,angle=0,scale=3.0]{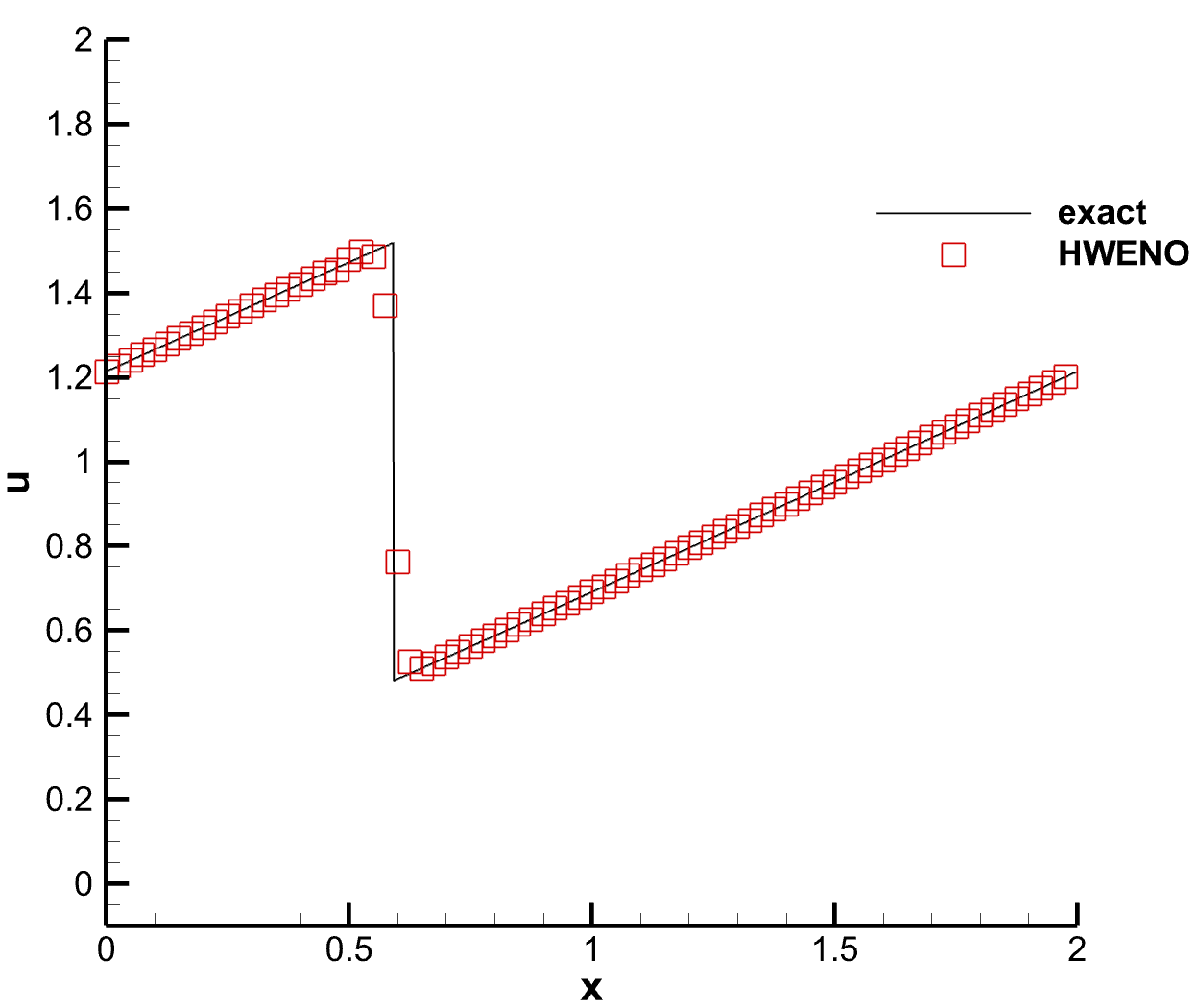}}
    \caption{Example 1: solutions with $N=80$, $C_a=0.01$.}
    \label{fig:1D Burgers'3}
\end{figure}

At \( T = 0.5/\pi \), the exact solution is smooth, and the corresponding $L^1$ errors and $L^{\infty}$ errors of $u$ are listed in Tables \ref{tab:1D Burgers'1} and \ref{tab:1D Burgers'2}, respectively. We can observe that our scheme is always stable and can achieve the designed fifth order for different $C_a$.
A discontinuity appears in the exact solution at \( T = 1.5/\pi \), and the discontinuity has propagated through the left boundary by \( T = 5/\pi \). The numerical solutions with $N=80$ and exact solutions at these times are compared in Figure \ref{fig:1D Burgers'3}, showing that the proposed schemes can capture the discontinuity well without noticeable oscillations.

\begin{table}[htb!]
	\centering
        \setlength{\tabcolsep}{4mm}
    \caption{Example 1: $L^1$ errors and orders of accuracy of $u$ at $T=0.5/\pi$.}
	\begin{tabular}{|c|cc|cc|cc|}\hline
        & \multicolumn{2}{c|}{$C_a=0.01$ }& \multicolumn{2}{c|}{$C_a=0.50$ }& \multicolumn{2}{c|}{$C_a=0.99$ } \\ \cline{2-7}
		N& $L^1$ error  & order & $L^1$ error &order& $L^1$ error &order\\ \hline
		10  &4.274E-02  & -     &2.558E-02  & -    & 4.099E-02 & -\\ 
		20  &3.543E-03  & 3.59 &6.995E-03  & 1.87& 2.976E-03 & 3.78\\ 
		40  &5.591E-05  & 5.99 &4.561E-04  & 3.94& 9.726E-05 & 4.94\\ 
		80  &1.838E-06  & 4.93 &1.515E-05  & 4.91& 1.949E-06 & 5.64\\ 
		160 &2.393E-08  & 6.26 &3.791E-08  & 8.64& 2.091E-08 & 6.54\\ 
		320 &5.928E-10  & 5.34 &5.119E-10  & 6.21& 5.054E-10 & 5.37\\ 
		640 &1.739E-11  & 5.09 &1.462E-11  & 5.13& 1.490E-11 & 5.08\\ \hline
	\end{tabular}
 \label{tab:1D Burgers'1}
\end{table}

\begin{table}[htb!]
	\centering
        \setlength{\tabcolsep}{4mm}
    \caption{Example 1: $L^\infty$ errors and orders of accuracy of $u$ at $T=0.5/\pi$. }
	\begin{tabular}{|c|cc|cc|cc|}\hline
        & \multicolumn{2}{c|}{$C_a=0.01$ }& \multicolumn{2}{c|}{$C_a=0.50$ }& \multicolumn{2}{c|}{$C_a=0.99$ } \\ \cline{2-7}
		N& $L^\infty$ error  & order & $L^\infty$ error &order& $L^\infty$ error &order\\ \hline
		10  &8.752E-02  & -     &8.223E-02  & -    & 8.421E-02 & -\\ 
		20  &1.771E-02  & 2.30 &1.859E-02  & 2.15& 1.470E-02 & 2.52\\ 
		40  &2.542E-04  & 6.12 &3.372E-03  & 2.46& 7.600E-04 & 4.27\\ 
		80  &2.369E-05  & 3.42 &1.687E-04  & 4.32& 3.704E-05 & 4.36\\ 
		160 &2.215E-07  & 6.74 &5.792E-07  & 8.19& 2.216E-07 & 7.38\\ 
		320 &6.338E-09  & 5.13 &6.340E-09  & 6.51& 6.339E-09  & 5.13\\ 
		640 &1.890E-10  & 5.07  &1.889E-10  & 5.07& 1.890E-10 & 5.07\\ \hline
	\end{tabular}
 \label{tab:1D Burgers'2}
\end{table}

\vspace{0.3cm}
\noindent
\textbf{Example 2.}
Here, we test the one-dimensional linear system 
\[
\begin{cases}
    u_t+w_x=0,\quad x\in(0,2\pi),\,\, t>0,\\
    w_t+u_x=0,\quad x\in(0,2\pi),\,\, t>0,\\
    u(x,0)=\sin(x),\quad
    w(x,0)=-\sin(x),\quad x\in[0,2\pi],\\
    u(0,t)=\sin(t),\quad
    u(2\pi,t)=\sin(t),\quad t>0.
\end{cases}
\]
We take $C_b = 1-C_a$ and \( \Delta x = 2 \pi / N \). In this example, only one boundary condition is required at both \( x = 0 \) and \( x = 2\pi \). Here, Dirichlet boundary conditions for \( u \) are imposed. The exact solution is \( u(x,t) = -w(x,t) = \sin( x + t ) \). 
We provide the \( L^1 \) error in Table \ref{tab:1D eqs1} and the \( L^\infty \) error in Table \ref{tab:1D eqs2} for \( u \) at \( T = 1 \) with respect to different values of \( C_a \), demonstrating that the scheme is stable and can reach the designed high order for linear system.

\begin{table}[htb!]
	\centering
        \setlength{\tabcolsep}{4mm}
    \caption{Example 2: $L^1$ errors and orders of accuracy of $u$ at $T=1$.}
	\begin{tabular}{|c|cc|cc|cc|}\hline
        & \multicolumn{2}{c|}{$C_a=0.01$ }& \multicolumn{2}{c|}{$C_a=0.5$ }& \multicolumn{2}{c|}{$C_a=0.99$ } \\ \cline{2-7}
		N& $L^1$ error  & order & $L^1$ error &order& $L^1$ error &order\\ \hline
		10  &4.609E-02  & -     &2.870E-02  & -    & 4.606E-002 & -\\ 
		20  &1.503E-03  & 4.94 &1.417E-03  & 4.34& 1.519E-03 & 4.92\\ 
		40  &5.957E-05  & 4.66 &4.604E-05  & 4.94& 2.040E-05 & 6.22\\ 
		80  &1.945E-06  & 4.94 &1.197E-06  & 5.27& 1.272E-06 & 4.00\\ 
		160 &4.465E-08  & 5.45 &9.314E-09  & 7.01& 2.298E-08 & 5.79\\ 
		320 &1.526E-09  & 4.87 &2.903E-10  & 5.00& 7.882E-10 & 4.87\\ 
		640 &5.115E-11  & 4.90 &9.772E-12  & 4.89& 2.701E-11 & 4.87\\ \hline
	\end{tabular}
 \label{tab:1D eqs1}
\end{table}

\begin{table}[htb!]
	\centering
        \setlength{\tabcolsep}{4mm}
    \caption{Example 2: $L^{\infty}$ errors and orders of accuracy of $u$ at $T=1$.}
	\begin{tabular}{|c|cc|cc|cc|}\hline
        & \multicolumn{2}{c|}{$C_a=0.01$ }& \multicolumn{2}{c|}{$C_a=0.5$ }& \multicolumn{2}{c|}{$C_a=0.99$ } \\ \cline{2-7}
		N& $L^\infty$ error  & order & $L^\infty$ error &order& $L^\infty$ error &order\\ \hline
		10  &9.573E-02  & -    &7.534E-02  & -    & 1.571E-01 & -\\ 
		20  &4.987E-03  & 4.26 &3.100E-03  & 4.60& 1.364E-02 & 3.53\\ 
		40  &2.775E-04  & 4.17 &1.955E-04  & 3.99& 1.416E-04 & 6.59\\ 
		80  &1.894E-05  & 3.87 &9.778E-06  & 4.32& 8.972E-06 & 3.98\\ 
		160 &3.060E-07  & 5.95 &8.292E-08  & 6.88& 1.482E-07 & 5.92\\ 
		320 &1.027E-08  & 4.90 &2.944E-09  & 4.82& 5.553E-09 & 4.74\\ 
		640 &3.384E-10  & 4.92 &1.009E-10  & 4.87& 1.969E-10 & 4.82\\ \hline
	\end{tabular}
 \label{tab:1D eqs2}
\end{table}

\vspace{0.3cm}
\noindent
\textbf{Example 3.}
Next, we consider a 1D Euler equation on $[0,2\pi]$ with initial data
\[
\begin{cases}
    \rho(x,0)=1+0.2\sin(x),\\
    p(x,0)=2,\\
    u(x,0)=1.
\end{cases}
\]
The grid partitioning is the same as in Example 2. 
Boundary conditions on \( \rho \) and \( u \) at \( x = 0 \) and \( \rho \) at \( x = 2\pi \) are given such that 
the system has an exact solution \( \rho(x,t) = 1 + 0.2 \sin(x - t), p(x,t) = 2, u(x,t) = 1 \). At \( T = 1 \), the errors and orders of accuracy of the density \( \rho \) are provided in Tables \ref{tab:1D Euler1} and \ref{tab:1D Euler2}, indicating that our methods can achieve the designed fifth order accuracy.

\begin{table}[htb!]
	\centering
        \setlength{\tabcolsep}{4mm}
    \caption{Example 3: $L^1$ errors and orders of accuracy of $\rho$ at $T=1$.}
	\begin{tabular}{|c|cc|cc|cc|}\hline
        & \multicolumn{2}{c|}{$C_a=0.01$ }& \multicolumn{2}{c|}{$C_a=0.5$ }& \multicolumn{2}{c|}{$C_a=0.99$ } \\ \cline{2-7}
		N& $L^1$ error  & order & $L^1$ error &order& $L^1$ error &order\\ \hline
		10  &3.488E-03  & -     &1.585E-03  & -  & 4.784E-03 & -\\ 
		20  &1.481E-04  & 4.56 &2.225E-04  & 2.83& 9.988E-05 & 5.58\\ 
		40  &2.998E-06  & 5.63 &1.042E-06  & 7.74& 2.725E-06 & 5.20\\ 
		80  &7.021E-08  & 5.42 &1.662E-08  & 5.97& 3.382E-08 & 6.33\\ 
		160 &2.143E-09  & 5.03 &5.104E-10  & 5.03& 1.268E-09 & 4.74\\ 
		320 &7.642E-11  & 4.81 &1.735E-11  & 4.88& 4.635E-11 & 4.77\\ 
		640 &2.795E-12  & 4.77 &7.551E-13  & 4.52& 1.842E-12 & 4.65\\ \hline
	\end{tabular}
 \label{tab:1D Euler1}
\end{table}

\begin{table}[htb!]
	\centering
    \setlength{\tabcolsep}{4mm}
    \caption{Example 3: $L^\infty$ errors and orders of accuracy of $\rho$ at $T=1$.}
	\begin{tabular}{|c|cc|cc|cc|}\hline
        & \multicolumn{2}{c|}{$C_a=0.01$ }& \multicolumn{2}{c|}{$C_a=0.5$ }& \multicolumn{2}{c|}{$C_a=0.99$ } \\ \cline{2-7}
		N& $L^\infty$ error  & order & $L^\infty$ error &order& $L^\infty$ error &order\\ \hline
		10  &1.417E-02  & -     &3.494E-03  & -    & 1.972E-02 & -\\ 
		20  &1.853E-03  & 2.93 &8.817E-04  & 1.99& 3.495E-04 & 5.82\\ 
		40  &6.883E-05  & 4.75 &4.081E-06  & 7.76& 2.177E-05 & 4.00\\ 
		80  &1.358E-06  & 5.66 &2.262E-07  & 4.32& 6.742E-07 & 5.01\\ 
		160 &3.630E-08  & 5.22 &7.956E-09  & 4.83& 2.606E-08 & 4.69\\ 
		320 &1.419E-09  & 4.68 &2.980E-10  & 4.74& 9.530E-10 & 4.77\\ 
		640 &5.128E-11  & 4.79 &1.044E-11  & 4.83& 3.251E-11 & 4.87\\ \hline
	\end{tabular}
 \label{tab:1D Euler2}
\end{table}

\begin{figure}[htb!]
    \centering
    \subfigure[$\Delta x$=1/400]
    {\includegraphics[width=1.0in,angle=0,scale=3.0]{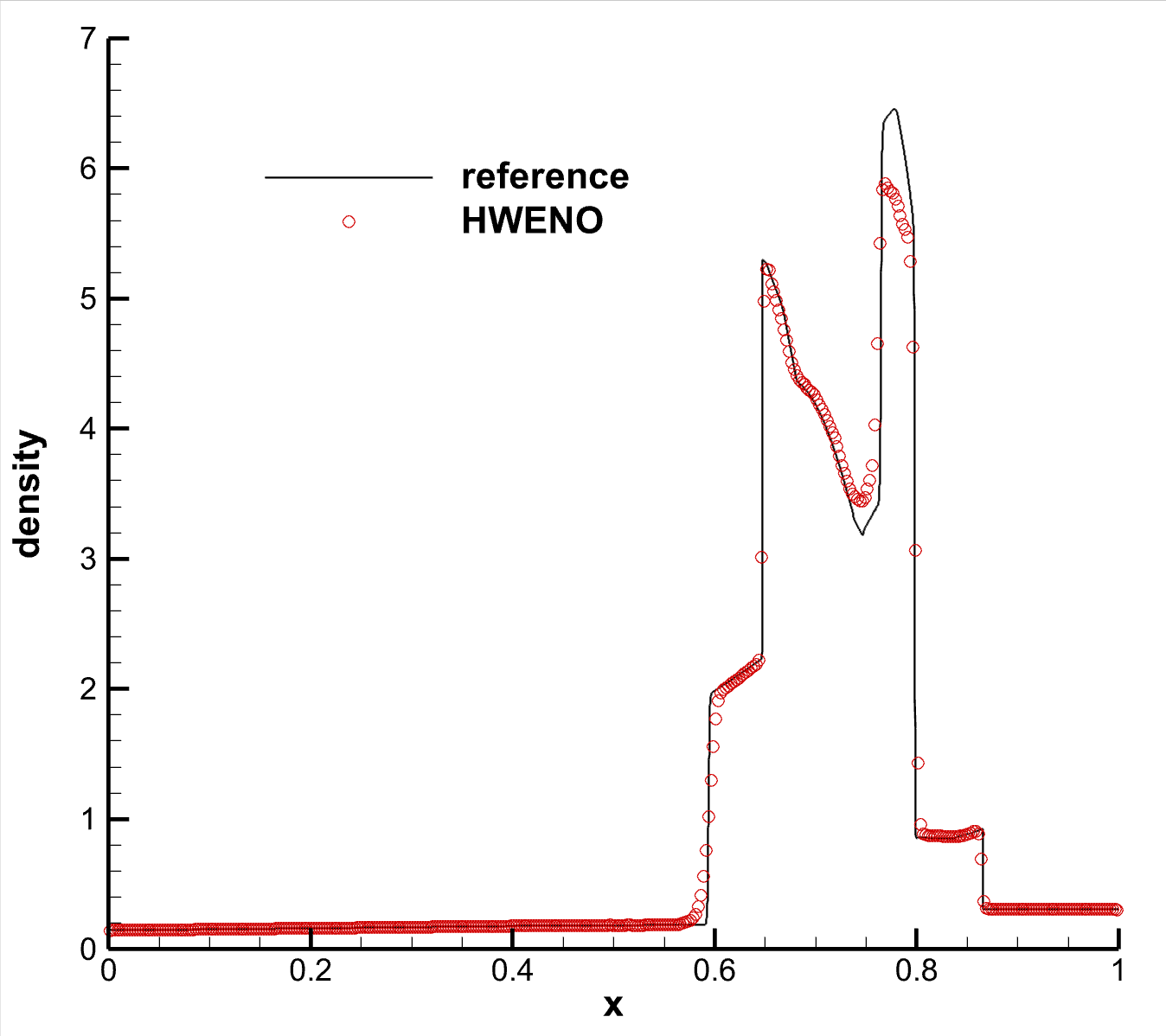}}
    \subfigure[$\Delta x$=1/800]
    {\includegraphics[width=1.0in,angle=0,scale=3.0]{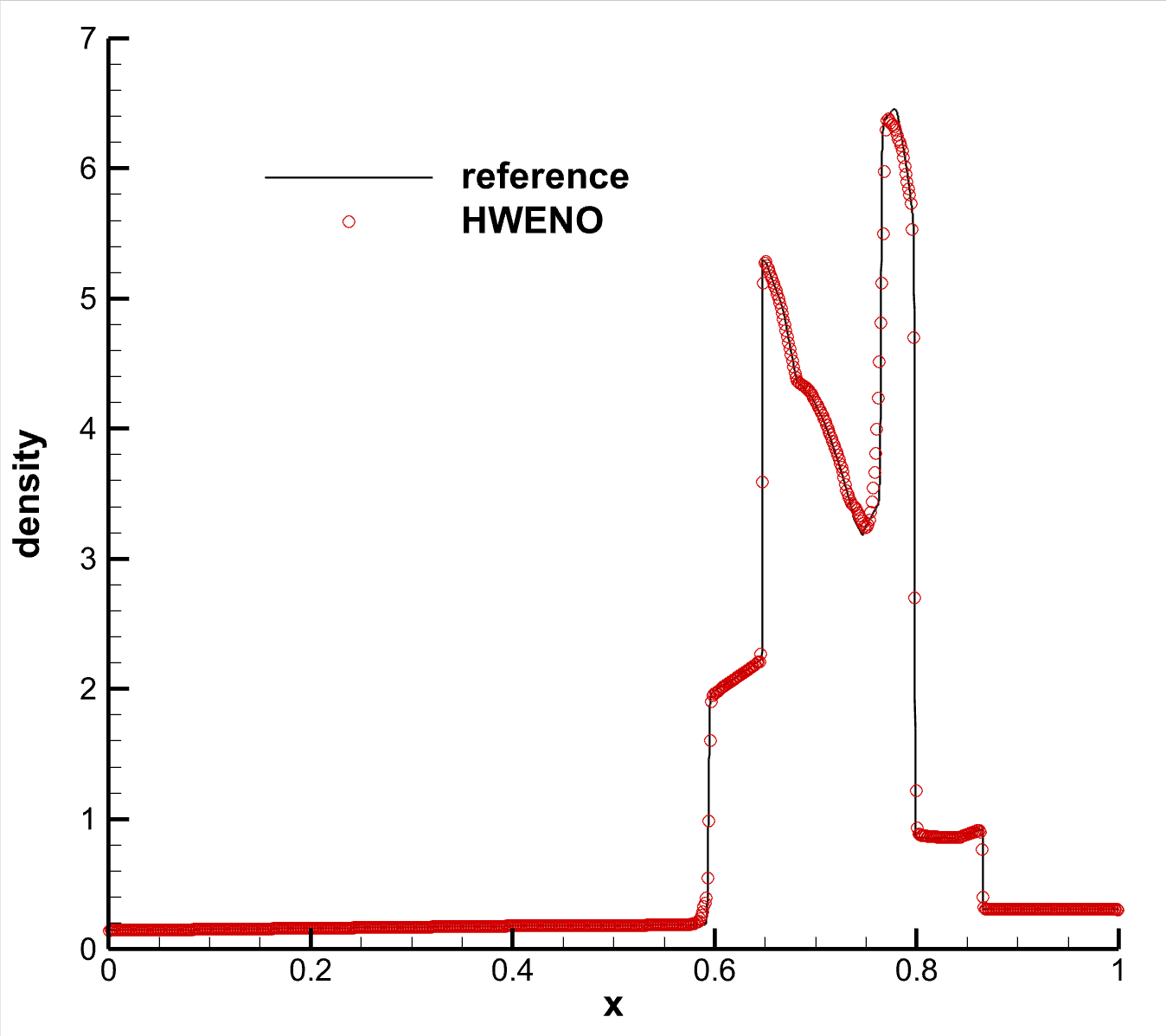}}
    \caption{Example 4: the density profiles of the blast wave problem.}
    \label{fig:blast}
\end{figure}

\vspace{0.3cm}
\noindent
\textbf{Example 4.}
Then we consider the interaction of two blast waves. In this problem, multiple reflections occur from shock waves and rarefaction waves on the wall. The initial conditions are
\[
(\rho,u,p)=\begin{cases}
    (1,0,1000),\quad&x\in[0,0.1),\\
    (1,0,0.01),\quad&x\in[0.1,0.9),\\
    (1,0,100),\quad&x\in[0.9,1.0].\\
\end{cases}
\]
We apply reflective boundary conditions at both boundaries and compute the numerical solution at \( T = 0.038 \). The reference solution is computed using the fifth order WENO scheme on body-fitted uniform mesh with \( \Delta x = 1/16000 \). 
Numerical results with different mesh size are provided in Figure \ref{fig:blast}. We can see that the proposed scheme can
distinguish the structure of the solution well.

\subsection{Two dimensional cases}
\begin{figure}[htb!]
    \centering
    \subfigure[Cut of the solution along $x=y$.]
    {\includegraphics[width=1.0in,angle=0,scale=3.0]{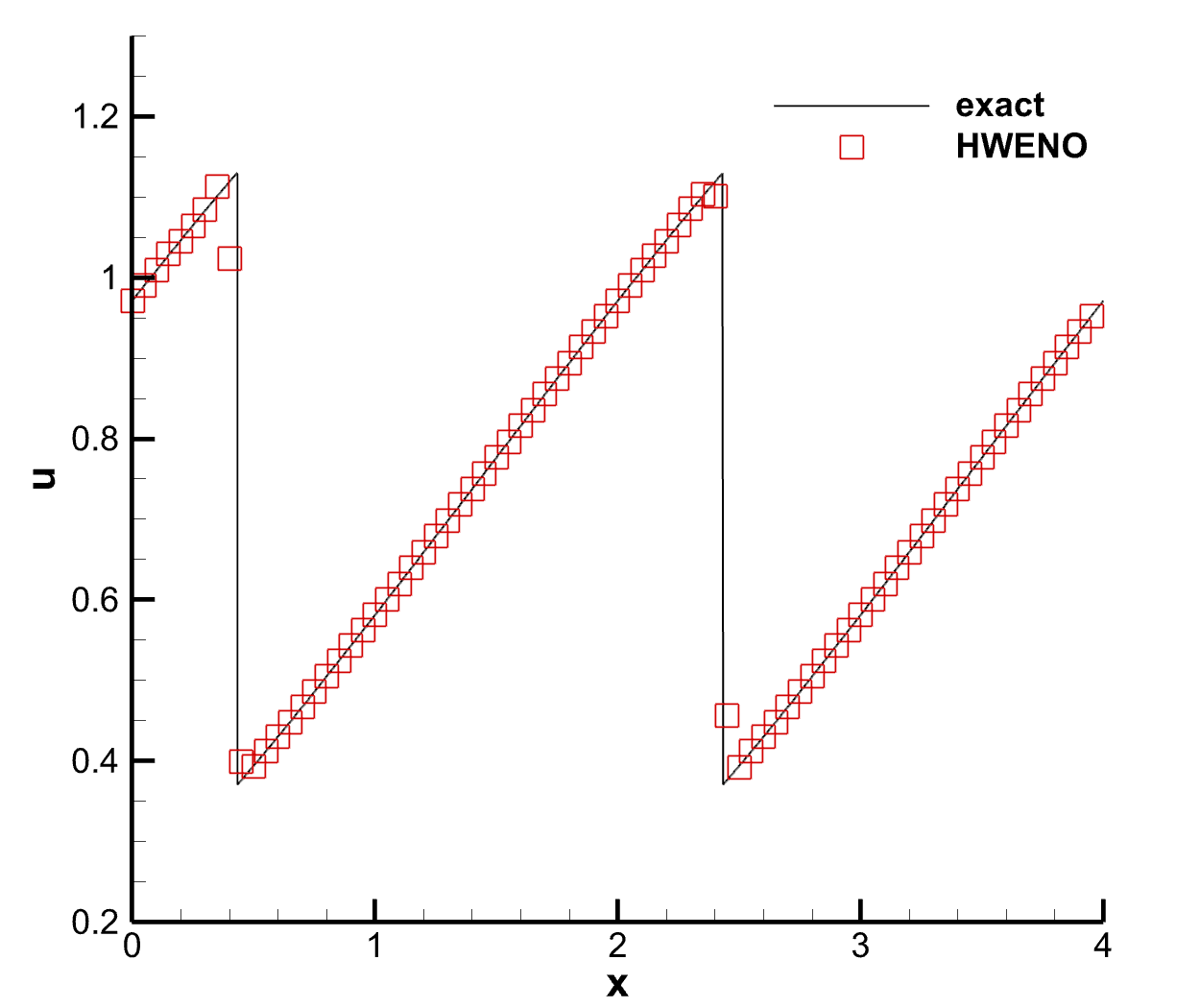}}
    \subfigure[Surface of the solution.]
    {\includegraphics[width=1.0in,angle=0,scale=3.0]{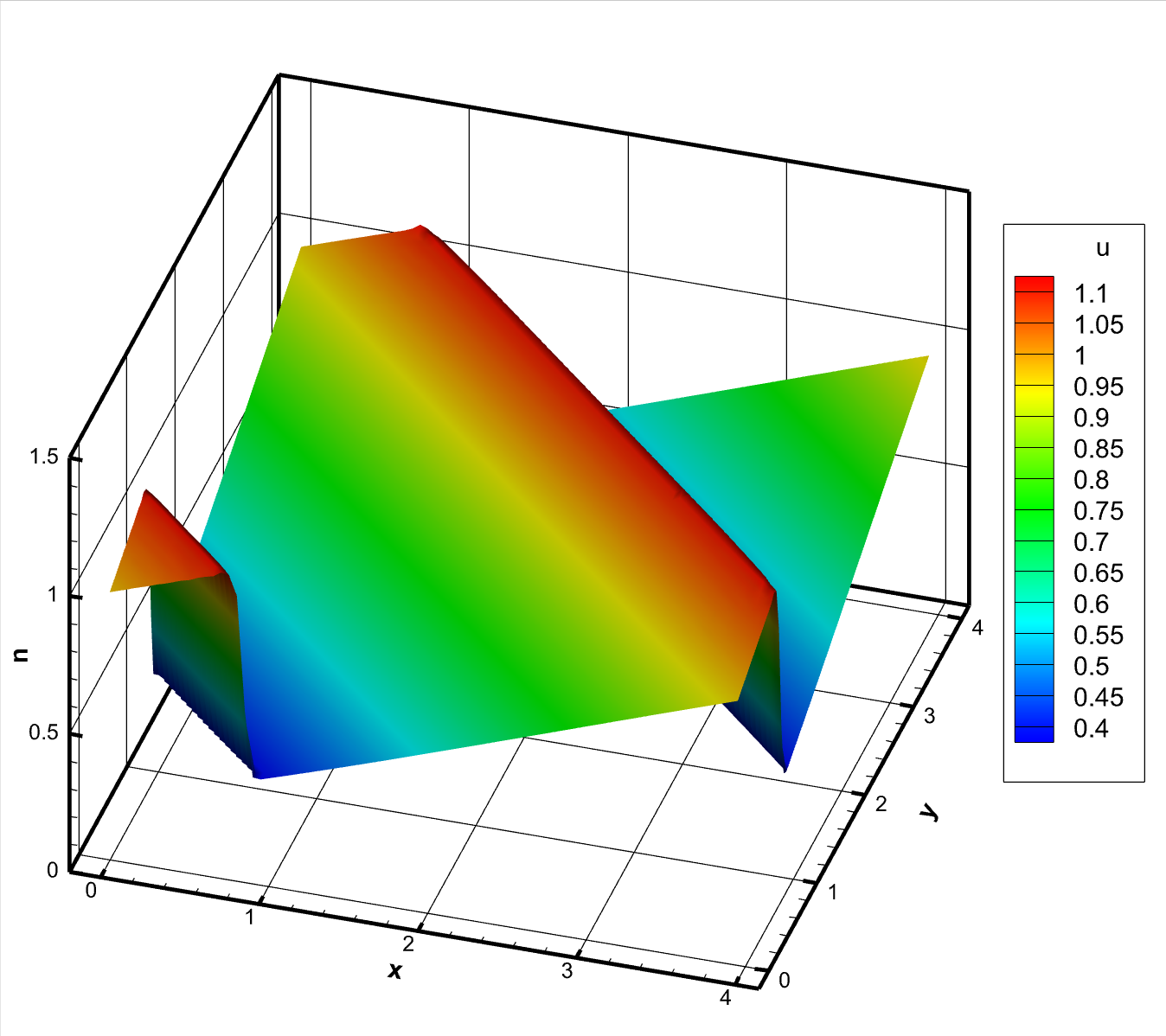}}
    \caption{Example 5: Solutions with $N_x\times N_y =80\times80$, $C_a=0.01$, $T=6/\pi$. }
    \label{fig:2D Burgers'3}
\end{figure}

\noindent
\textbf{Example 5.}
At first, we look at the 2D Burgers' equation on a rectangular domain
\[
\begin{cases}
    u_t+\left(\frac{u^2}{2}\right)_x+\left(\frac{u^2}{2}\right)_y=0,\quad (x,y)\in[0,4]\times[0,4],\quad t>0,\\
    u(x,y,0)=0.75+0.5\sin(\pi\left(\frac{x+y}{2}\right)),\quad (x,y)\in[0,4]\times[0,4].
\end{cases}
\]
Here we define \( x_i = (C_a - 1 + i) \Delta x\), \(y_i = (C_a - 1 + i) \Delta y\), \( 1 \leq i \leq N \), where \( C_a \in [0, 1) \) and \( \Delta x = \Delta y = 4/N \). 
The boundaries \( x = 0 \) and \( y = 0 \) are inflow boundaries, where the boundary conditions obtained from the exact solution with periodic boundaries are imposed here. On the other hand, \( x = 4 \) and \( y = 4 \) are outflow boundaries, which require no boundary conditions. 
The solution remains continuous at \( T = 1/\pi \), and the corresponding errors are listed in Tables \ref{tab:2D Burgers'1} and \ref{tab:2D Burgers'2}. 
Both the 1D cross-section along the diagonal and the two-dimensional plane at \( T = 6/\pi \) are presented in the figure \ref{fig:2D Burgers'3}, when discontinuities exist. 
We can see that our 2D method can achieve the optimal fifth order for smooth solutions and avoid numerical oscillation near discontinuities.

\begin{table}[htb!]
	\centering
    \setlength{\tabcolsep}{4mm}
    \caption{Example 5: $L^1$ errors and orders of accuracy of $u$ at $T=1/\pi$. }
	\begin{tabular}{|c|cc|cc|cc|}\hline
        & \multicolumn{2}{c|}{$C_a=0.01$ }& \multicolumn{2}{c|}{$C_a=0.5$ }& \multicolumn{2}{c|}{$C_a=0.99$ } \\ \cline{2-7}
		$N_x\times N_y$ & $L^1$ error  & order & $L^1$ error &order& $L^1$ error &order\\ \hline
		10$\times$10  &1.113E-02  & -     &1.735E-02 & -   & 6.862E-03 & -\\ 
		20$\times$20  &8.105E-04  & 3.78 &6.016E-04 & 4.85& 5.454E-04 & 3.65\\ 
		40$\times$40  &3.367E-05  & 4.59 &9.622E-06 & 5.97& 2.836E-05 & 4.27\\ 
		80$\times$80  &6.919E-07  & 5.60 &6.507E-07 & 3.89& 4.063E-07 & 6.13\\ 
		160$\times$160&1.114E-08  & 5.96 &1.067E-08 & 5.93& 1.017E-08 & 5.32\\ 
		320$\times$320&3.158E-10  & 5.14 &3.101E-10 & 5.10& 2.919E-10 & 5.12\\
        640$\times$640&9.349E-12  & 5.08 &9.179E-12 & 5.08& 8.858E-12 & 5.04\\\hline
	\end{tabular}
 \label{tab:2D Burgers'1}
\end{table}

\begin{table}[htb!]
	\centering
    \setlength{\tabcolsep}{4mm}
    \caption{Example 5: $L^\infty$ errors and orders of accuracy of $u$ at $T=1/\pi$}
	\begin{tabular}{|c|cc|cc|cc|}\hline
        & \multicolumn{2}{c|}{$C_a=0.01$ }& \multicolumn{2}{c|}{$C_a=0.5$ }& \multicolumn{2}{c|}{$C_a=0.99$ } \\ \cline{2-7}
		$N_x\times N_y$& $L^\infty$ error  & order & $L^\infty$ error &order& $L^\infty$ error &order\\ \hline
		10$\times$10  &4.094E-02  & -    &3.722E-02  & -   & 1.983E-02 & -\\ 
		20$\times$20  &2.896E-03  & 3.82 &2.449E-03  & 3.93& 2.852E-03 & 2.80\\ 
		40$\times$40  &1.548E-04  & 4.23 &8.994E-05  & 4.77& 1.528E-04 & 4.22\\ 
		80$\times$80  &1.003E-05  & 3.95 &7.778E-06  & 3.53& 4.573E-06 & 5.06\\ 
		160$\times$160&1.311E-07  & 6.26 &1.281E-07  & 5.92& 1.307E-07 & 5.13\\ 
		320$\times$320&3.857E-09  & 5.09 &3.876E-09  & 5.05& 3.852E-09 & 5.08\\
        640$\times$640&1.161E-10  & 5.05 &1.168E-10 & 5.05& 1.161E-10 & 5.05\\\hline
	\end{tabular}
 \label{tab:2D Burgers'2}
\end{table}

\begin{figure}[htb!]
    \centering
    \subfigure[Cut of the solution along $x=y$.]
    {\includegraphics[width=1.0in,angle=0,scale=3.0]{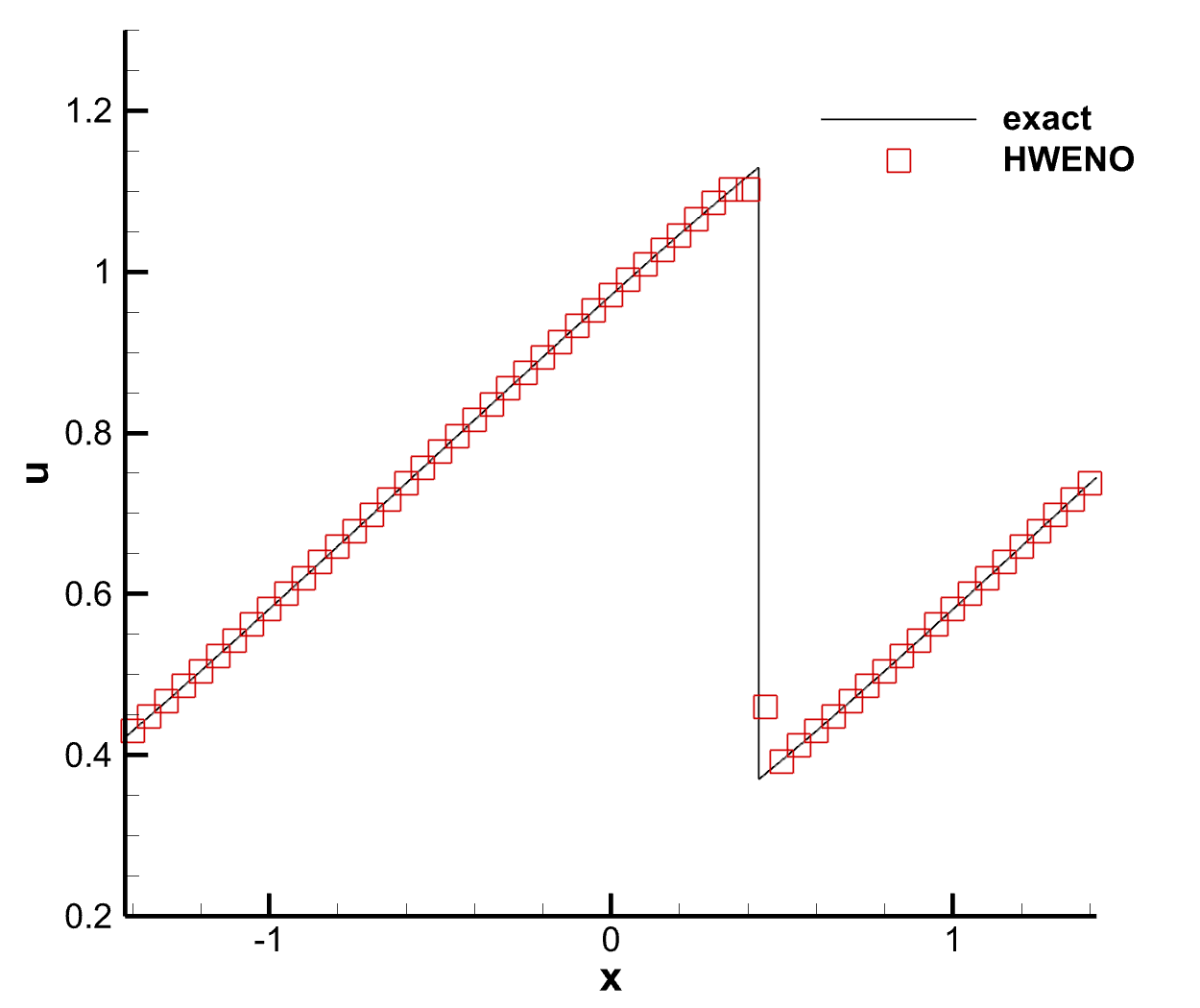}}
    \subfigure[Surface of the solution.]
    {\includegraphics[width=1.0in,angle=0,scale=3.0]{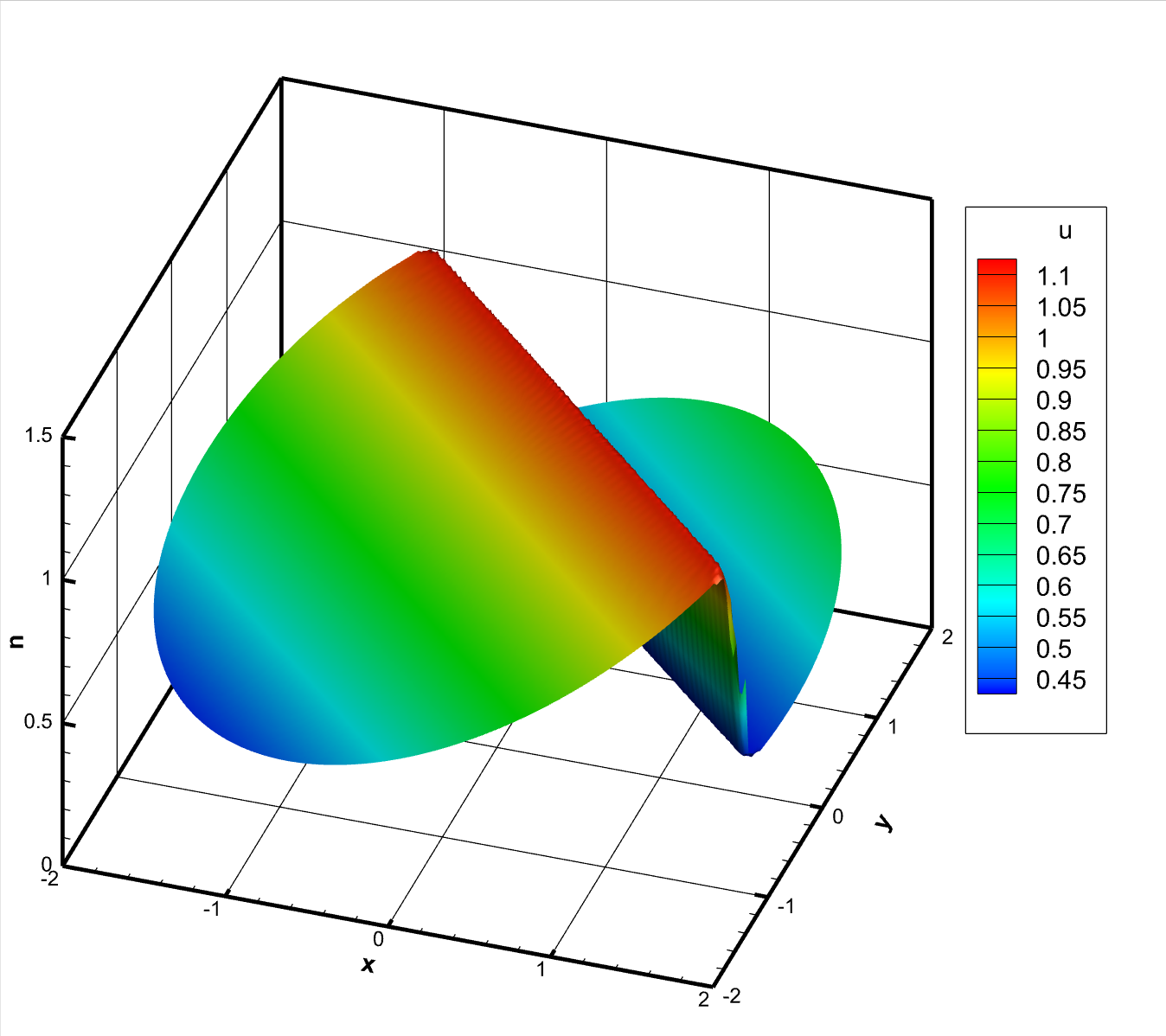}}
    \caption{Example 6: solutions with $\Delta x=\Delta y = 1/20$ at $T=6/\pi$. }
    \label{fig:2D Burgers'5}
\end{figure}

\vspace{0.3cm}
\noindent
\textbf{Example 6.}
Next, we consider the 2D Burgers' equation on a disk,
\[
\begin{cases}
    u_t+\left(\frac{u^2}{2}\right)_x+\left(\frac{u^2}{2}\right)_y=0,\quad (x,y)\in\Omega,\quad t>0,\\
    u(x,y,0)=0.75+0.5\sin(\pi\left(\frac{x+y}{2}\right)),\quad (x,y)\in\Omega,
\end{cases}
\]
with 
$$\Omega = \{(x, y) \mid x^2 + y^2 \leq 4\}.$$ 
Here, we define \( x_i = i \Delta x \) and \( y_i = i \Delta y \). We impose the problem with Dirichlet boundary conditions on the inflow boundary such that the exact solution is the same as that of Example 5. 
The errors of the numerical solution at \( t = 1/\pi \) are provided in Table \ref{tab:2D Burgers'4}, and the plots along the diagonal as well as the 2D plane are shown in Figure \ref{fig:2D Burgers'5}. 
It is observed that our scheme is effective for curved boundary. 

\begin{table}[H]
	\centering
        \setlength{\tabcolsep}{4mm}
	\caption{Example 6: errors and orders of accuracy of $u$ at $T=1/\pi$.}
	\begin{tabular}{|c|cccc|}\hline
		$\Delta x = \Delta y$ &  $L^1$ error  & order &  $L^\infty$ error  & order\\ \hline
		1/5  &1.173E-03  & - &4.202E-03  & -\\ 
		1/10  &5.161E-05  & 4.51 &2.727E-04  & 3.95\\ 
		1/20  &1.347E-06  & 5.26 &2.129E-05  & 3.68\\ 
		1/40 &1.571E-08  & 6.42 &1.309E-07  & 7.35\\ 
		1/80 &4.514E-10  & 5.12 &8.986E-09  & 3.86\\
        1/160 &1.284E-11  & 5.14 &2.444E-10  & 5.20\\\hline
	\end{tabular}
 \label{tab:2D Burgers'4}
\end{table}

\begin{figure}[htb!]
    \centering
    \subfigure[ Pressure contour.]
    {\includegraphics[width=1.2in,angle=0,scale=2.0]{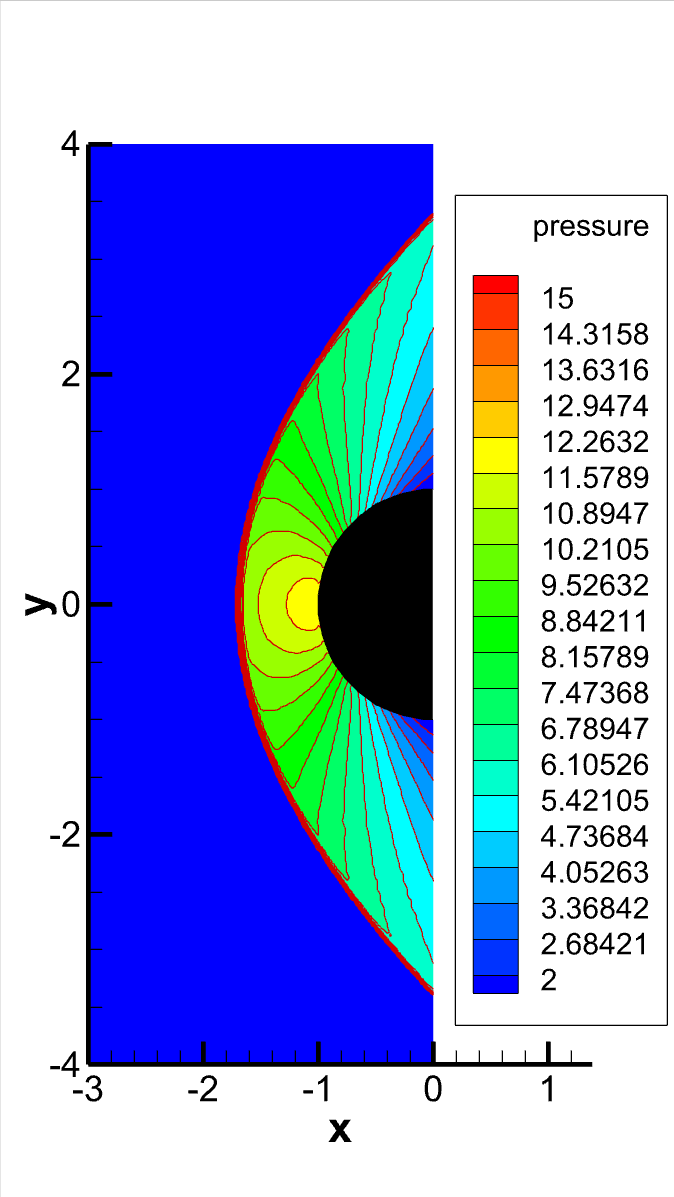}}
    \subfigure[Mach number contour.]
    {\includegraphics[width=1.2in,angle=0,scale=2.0]{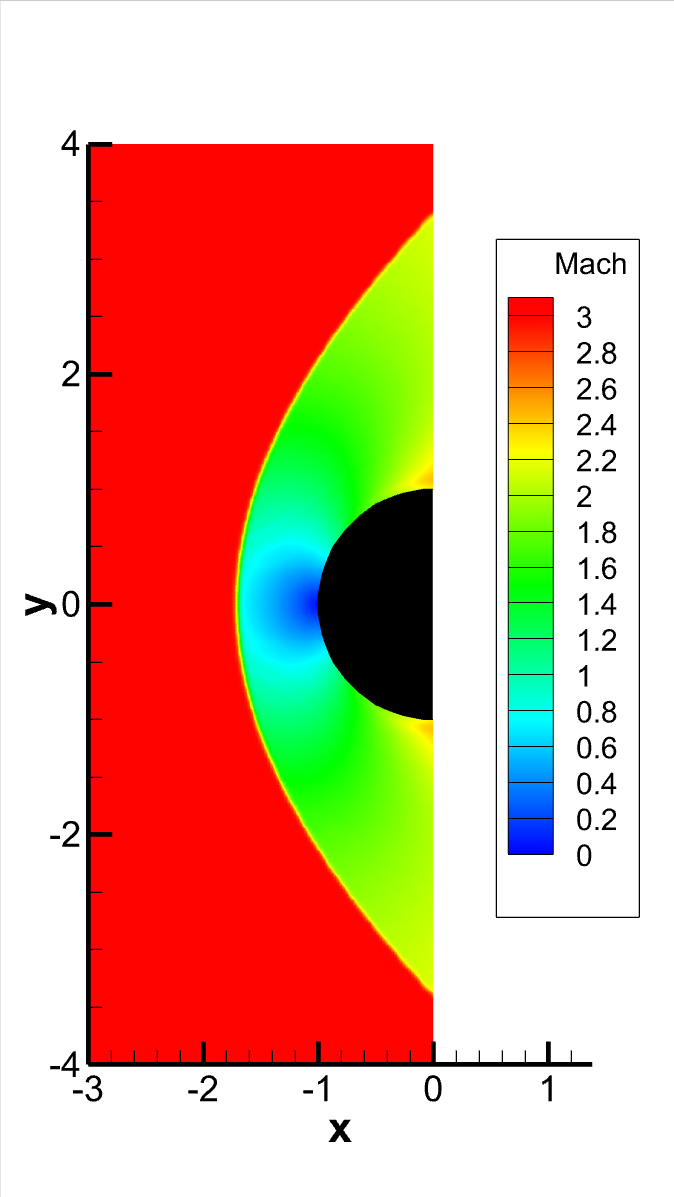}}
    \caption{Example 7: pressure and Mach number of flow past a cylinder with $\Delta x=\Delta y=1/40$.}
    \label{fig:cylinder}
\end{figure}

\vspace{0.3cm}
\noindent
\textbf{Example 7.}
We use an example with a curved wall to test our schemes. The initialization of the problem is a Mach 3 flow moving from the left to a unit radius cylinder located at the origin of the $x-y$ plane. Due to the symmetry of the problem, we only need to compute the upper half plane, i.e., $[-3,0]\times [0,6]$. On the surface of the cylinder, we use the ILW method to deal with the solid-wall boundary conditions. At the left boundary of the computed region $x=-3$, the inflow boundary conditions are given. And for the right boundary $x=0$ and the upper boundary $y=6$, outflow boundary conditions are imposed.
Pressure and Mach number are given in Figure \ref{fig:cylinder}, which are comparable with those in the previous work. 

\begin{figure}[htb!]
	\centering
	\includegraphics[width=0.6\linewidth]{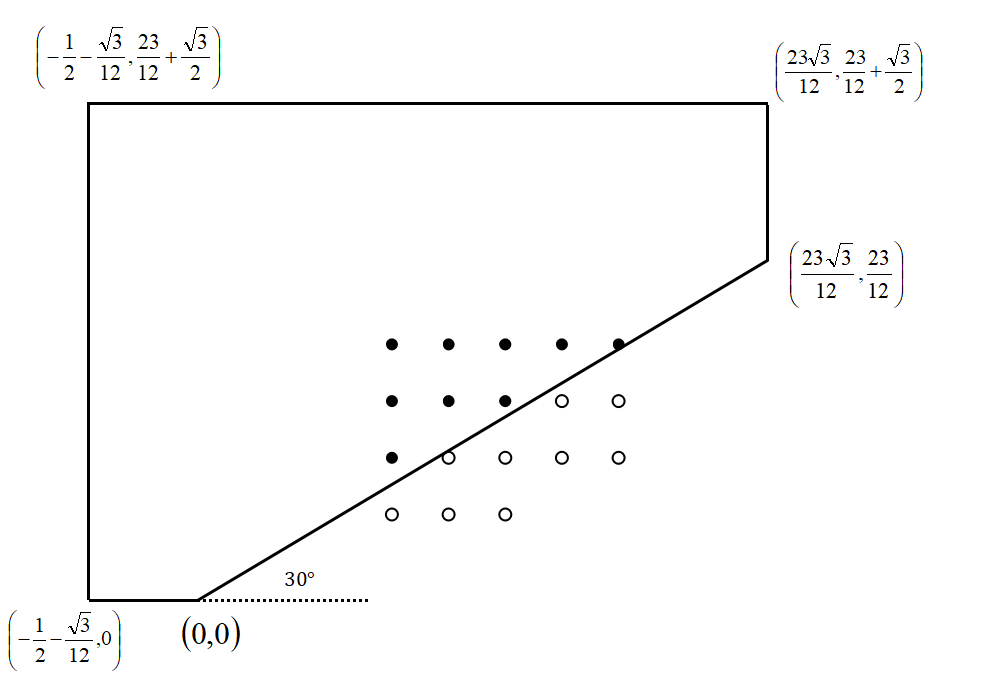}
	\caption{Example 8: double Mach reflection computational region}
	\label{fig:DMR}
\end{figure}

\begin{figure}[htb!]
    \centering
    \subfigure[$\Delta x=\Delta y=1/320$]{\includegraphics[width=1.0in,angle=0,scale=3.0]{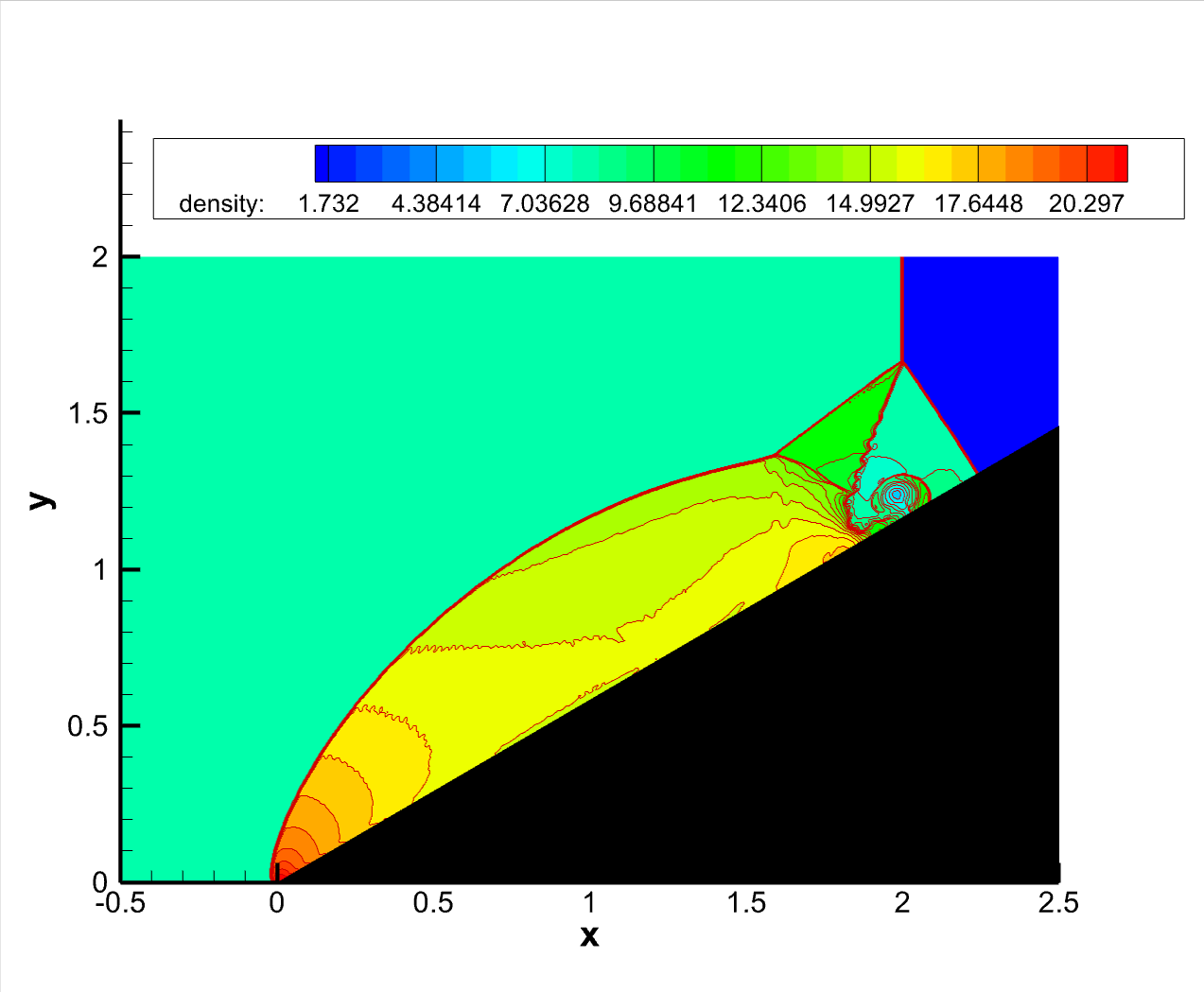}}
    \subfigure[$\Delta x=\Delta y=1/320$, zoomed-in figure.]{\includegraphics[width=1.0in,angle=0,scale=3.0]{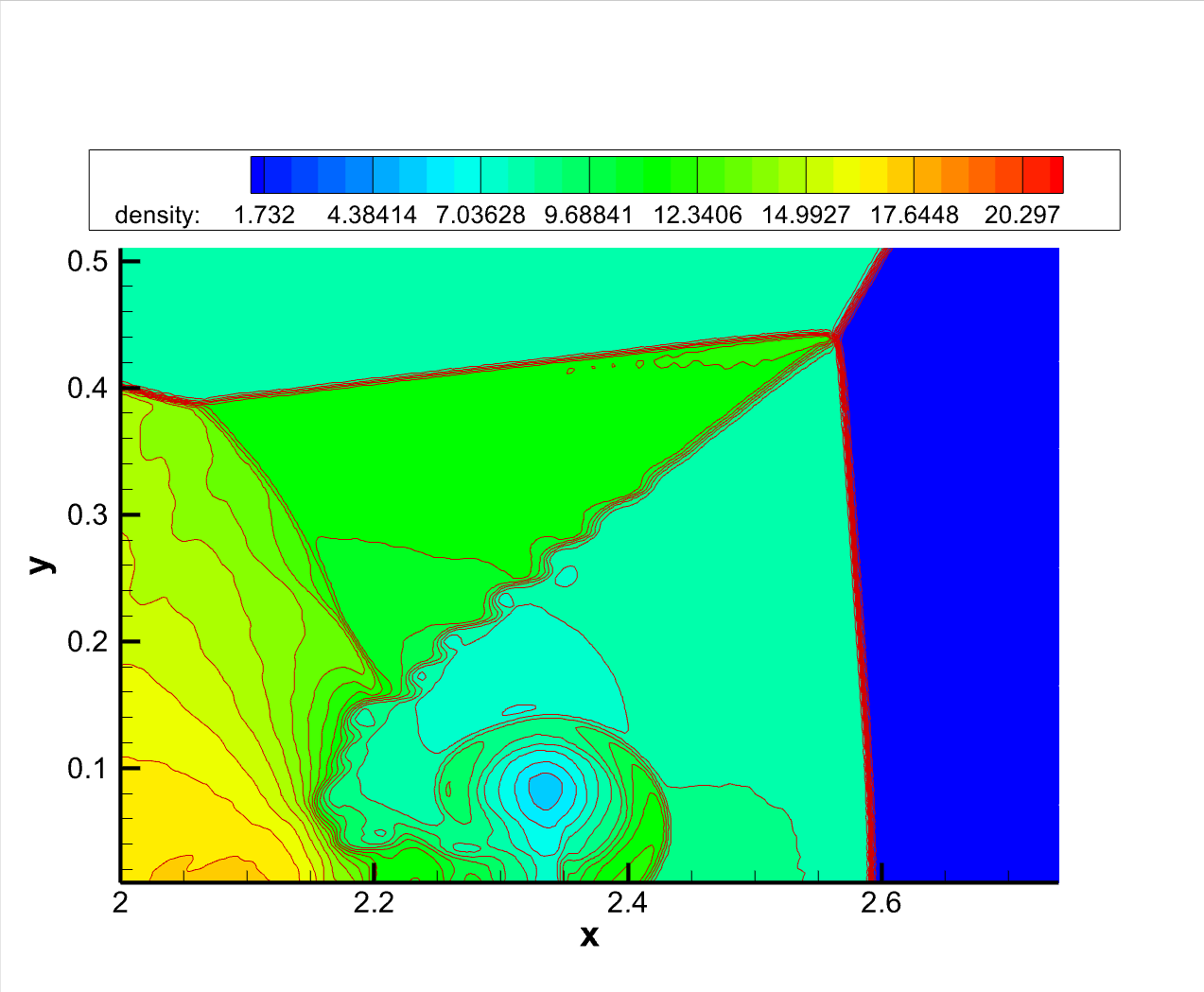}}
    \subfigure[$\Delta x=\Delta y=1/640$]{\includegraphics[width=1.0in,angle=0,scale=3.0]{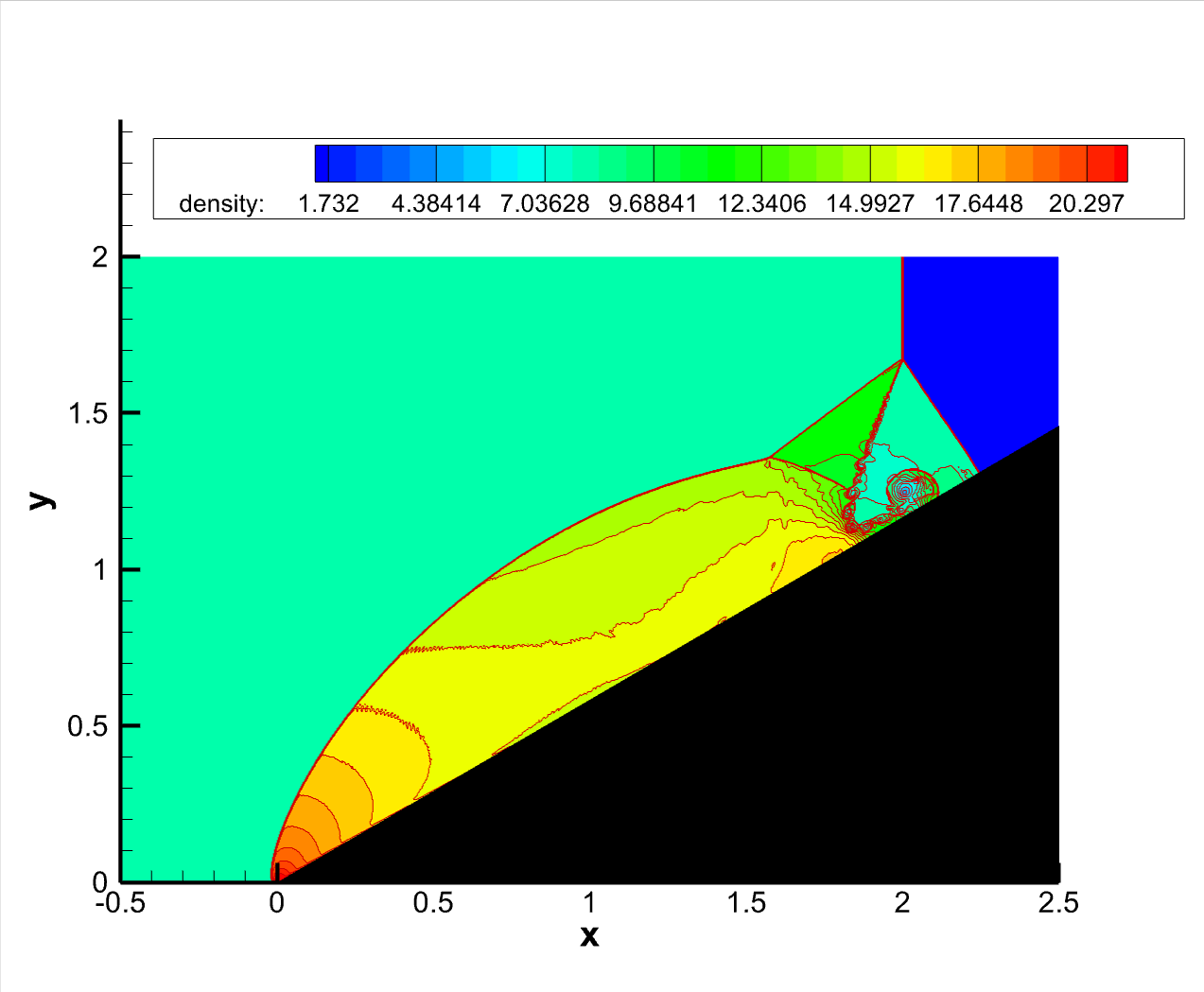}}
    \subfigure[$\Delta x=\Delta y=1/640$, zoomed-in figure.]{\includegraphics[width=1.0in,angle=0,scale=3.0]{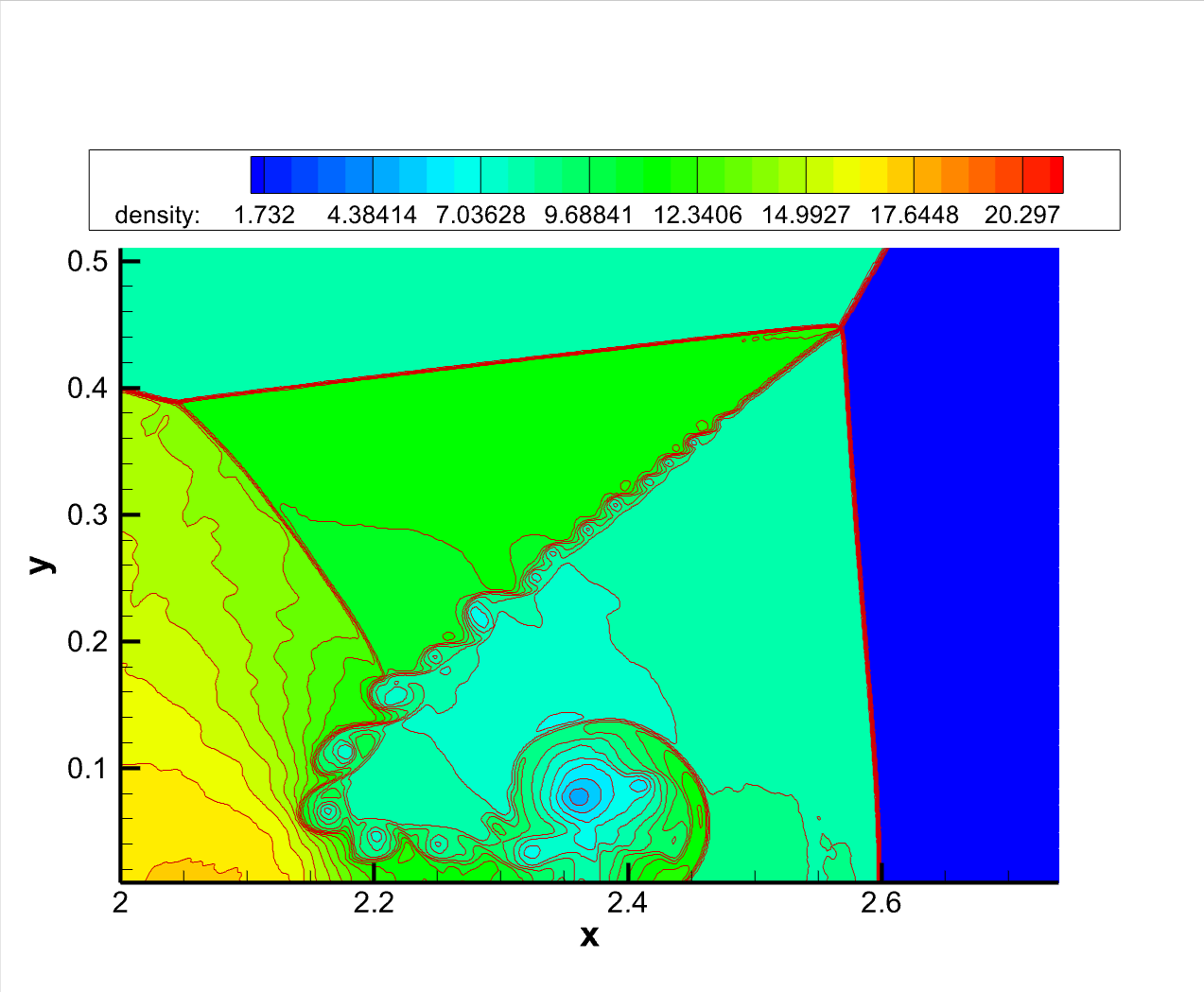}}
    \caption{Example 8: density contour of double Mach reflection, 30 contours from 1.731 to 20.92. }
    \label{fig:ex6-2}
\end{figure}

\vspace{0.3cm}
\noindent
\textbf{Example 8.}
Finally, we test the double Mach reflection problem. This example describes the interaction between a shock wave with a Mach number of 10 moving in the horizontal direction and a solid wall at an angle of $30^\circ$. The domain is shown in Figure \ref{fig:DMR}. At the top boundary of the computational domain, we set appropriate boundary conditions based on the given Mach number, while supersonic inflow and outflow boundary conditions are applied on the left and right sides, respectively. For the hypotenuse solid wall, the designed ILW method is employed.
Density contours with different mesh sizes are presented in Figure \ref{fig:ex6-2}, and the zoomed-in regions near the double Mach stem are also given. It is observed that the results are comparable with the previous results, indicating the effectiveness of our ILW method. 

\section{Conclusion}

In this paper, we introduce a novel ILW boundary treatment method for fifth-order finite difference HWENO schemes \cite{fan2023robust} to solve hyperbolic conservation laws. 
Distinct from existing ILW approaches \cite{tan2010inverse, tan2012efficient, liu2024113259}, our method incorporates additional parameters \( k \) and employs a least squares method. By adjusting the parameter \( k \), our method could achieves stability while reducing the required terms with ILW procedure. Therefore our scheme enhances computational efficiency, particularly for multi-dimensional systems. 
Additionally, we provide a stability analysis for the fully discrete schemes, specifying appropriate parameter values to ensure the boundary treatment does not effect the original CFL condition constraint of the inner scheme, thus it can avoid the cut-cell problem effectively.  
Through numerical experiments, we demonstrate that our algorithm is stable, robust, and achieves the designed fifth-order accuracy, performing exceptionally well in complex problems such as cylindrical flow and double Mach reflection problems.


\appendix

\section{Fifth order HWENO algorithm} \label{sec:app1}

Here, we give the formulation of the fifth order HWENO scheme to construct $\hat{f}^{+}_{i+1/2}$ and $\hat{h}^{+}_{i+1/2}$ based on $\{(f^+_{j}, h^+_{j}), j= i-1, i, i+1\}$\cite{fan2023robust}. 

\begin{enumerate}
    \item Construct a quartic polynomial $p_0(x)$ and two quadratic polynomials $p_1(x),p_2(x)$, satisfying
\begin{equation*}
    \left\{
    \begin{aligned}
        &\frac{1}{\Delta x}\int_{I_{i+l}}p_0(x)dx=f^+_{i+l},\quad l=-1,0,1,\\
		&\frac{1}{\Delta x}\int_{I_{i+l}}p'_0(x)dx=h^+_{i+l},\quad l=-1,1.
    \end{aligned}
    \right.
\end{equation*}
\begin{equation*}
		\frac{1}{\Delta x}\int_{I_{i+l}}p_1(x)dx=f^+_{i+l},\quad l=-1,0,
\end{equation*}
\begin{equation*}
		\frac{1}{\Delta x}\int_{I_{i+l}}p_2(x)dx=f^+_{i+l},\quad l=0,1.
\end{equation*}
Compute the values and derivatives at half point $x_{i+\frac{1}{2}}$:
\begin{equation*}
    \begin{aligned}
        &p_0(x_{i+\frac{1}{2}})=-\frac{23}{120}f^+_{i-1}+\frac{19}{30}f^+_{i}+\frac{67}{120}f^+_{i+1}-\frac{3\Delta x}{40}h^+_{i-1}-\frac{7\Delta x}{40}h^+_{i+1},\\
        &p_1(x_{i+\frac{1}{2}})=-\frac{1}{2}f^+_{i-1}+\frac{3}{2}f^+_{i},\\
        &p_2(x_{i+\frac{1}{2}})=\frac{1}{2}f^+_{i}+\frac{1}{2}f^+_{i+1},\\
        &p'_0(x_{i+\frac{1}{2}})=\frac{1}{\Delta x}(\frac{3}{8}f^+_{i-1}-2f^+_{i}+\frac{13}{8}f^+_{i+1})+\frac{1}{8}h^+_{i-1}-\frac{3}{8}h^+_{i+1}.
    \end{aligned}
\end{equation*}
Define the linear weights $\gamma_l$, $l=0,1,2$, satisfying 
\[
\gamma_l\geq 0, \quad \text{and} \quad 
\gamma_0+\gamma_1+\gamma_2=1.
\]
Here, we take \(\gamma_0=0.99,\gamma_1=\gamma_2=0.005\).

\item Compute the smoothness indicator $\beta_l$ at interval $I_i$ to measure the smoothness of the polynomial $p_l(x)$, $l=0,1,2,$
\begin{equation*}
\beta_l=\sum_{\alpha=1}^{k}\int_{I_{i}}\Delta x^{2\alpha-1}\left( \frac{d^{\alpha}p_l(x)}{dx^{\alpha}}\right)^2dx,\quad l=0,1,2.
\end{equation*}

\item Form the nonlinear weight $\omega_l$ based on the linear weight $\gamma_l$ and the smooth indicator $\beta_l$:
\[
\omega_l=\frac{\bar{\omega}_l}{\sum_{l=0}^2\bar{\omega}_l},\quad
\bar{\omega}_l=\gamma_l\left( 1+\frac{\tau}{\beta_l+\epsilon} \right),\quad l=0,1,2
\]
with
\[
\tau=\frac{1}{4}(|\beta_0-\beta_1|+|\beta_0-\beta_2|)^2.
\]
And $\epsilon=10^{-10}>0$ is a tiny number to avoid the denominator to become 0.

\item Finally, we reconstruct $\hat{f}^+_{i+\frac{1}{2}}$ with nonlinear weights and $\hat{h}^+_{i+\frac{1}{2}}$ with linear weights:
\begin{align*}
    \hat{f}^+_{i+\frac{1}{2}}=& \omega_0 \left( \frac{1}{\gamma_0}p_0(x_{i+\frac{1}{2}})
    -\frac{\gamma_1}{\gamma_0}p_1(x_{i+\frac{1}{2}})
    -\frac{\gamma_2}{\gamma_0}p_2(x_{i+\frac{1}{2}}) \right)
    +\omega_1p_1(x_{i+\frac{1}{2}})
    +\omega_2p_2(x_{i+\frac{1}{2}}),\\
    \hat{h}^+_{i+\frac{1}{2}}= & 
    p_0'(x_{i+\frac{1}{2}})
\end{align*}
\end{enumerate}

\section{Modification of the derivative} \label{sec:app2}
Here, we present the modification on the derivative $v_i$ \cite{fan2023robust} in the following.

\begin{enumerate}
    
\item Firstly, construct a quartic polynomial $q_0(x)$ and two first-order polynomial $q_1(x),q_2(x)$:
\[
    \begin{cases}
        q_0(x_{i+l})=u_{i+l},\quad l=-1,0,1,\\
        q_0'(x_{i+l})=v_{i+l},\quad l=-1,1.
    \end{cases}
\]
\[q_1(x_{i+l})=u_{i+l},\quad l=-1,0,\]
\[ q_2(x_{i+l})=u_{i+l},\quad l=0,1.\]
Then we have
\begin{equation*}
    \begin{aligned}
        &q_0'(x_i)=\frac{3}{4\Delta x}(u_{i+1}-u_{i-1})-\frac{1}{4}(v_{i-1}+v_{i+1}),\\
        &q_1'(x_i)=\frac{1}{\Delta x}(u_i-u_{i-1}),\\
        &q_2'(x_i)=\frac{1}{\Delta x}(u_{i+1}-u_{i}).
    \end{aligned}
\end{equation*}

\item Compute the smooth indicator $\beta_l$ of the polynomial $q_0(x),q_1(x),q_2(x)$.
\begin{equation*}
\beta_l=\sum_{\alpha=1}^{k}\int_{I_{i}}\Delta x^{2\alpha-1}\left( \frac{d^{\alpha}q_l(x)}{dx^{\alpha}} \right)^2dx,\quad l=0,1,2.
\end{equation*}

\item Define the nonliear weights 
\[
\lambda_l=\frac{\bar{\lambda}_l}{\sum_{l=0}^2\bar{\lambda}_l},\quad
\bar{\lambda}_l=d_l\left(1+\frac{\tau}{\beta_l+\epsilon}\right),\quad l=0,1,2
\]
with linear weights \(d_0=0.9,d_1=d_2=0.05\) and 
\[
\tau=\frac{1}{4}(|\beta_0-\beta_1|+|\beta_0-\beta_2|)^2.
\]

\item Finally, obtain the modification of the first derivative at $x_{i}$ as
\[
\tilde{v}_{i}=\lambda_0 \left( \frac{1}{d_0}q_0'(x_{i})-\frac{d_1}{d_0}q_1'(x_{i} ) -\frac{d_2}{d_0}q_2'(x_{i}) \right) 
+\lambda_1q_1'(x_{i}) +\lambda_2q_2'(x_{i}).
\]
\end{enumerate}

\bibliographystyle{abbrv}
\bibliography{reference}

\end{document}